\documentclass{gtart}

\def\ifplaintex{\expandafter\ifx\csname documentclass\endcsname\relax}

\def\gtp{{\mathsurround=0pt\it $\cal G\mskip-2mu$eometry \&\ 
$\cal T\!\!$opology $\cal P\!$ublications}}  

\def\recd{{\small Received:\qua\receiveddate\ifx\reviseddate\relax
\else\qquad Revised:\qua\reviseddate\fi\par}} 

\def\lognumber#1{\def\thelognumber{#1}}
\def\volumenumber#1{\def\thevolumenumber{#1}}
\def\volumeyear#1{\def\thevolumeyear{#1}}
\def\papernumber#1{\def\thepapernumber{#1}}
\def\pagenumbers#1#2{\def\startpage{#1}\def\finishpage{#2}}
\def\published#1{\def\publishdate{#1}}

\def\received#1{\def\receiveddate{#1}}
\def\revised#1{\def\reviseddate{#1}}
\def\accepted#1{\def\accepteddate{#1}}
\def\asciititle#1{\def\theasciititle{#1}}

\def\asciiaddress#1{\def\theasciiaddress{#1}}

\long\def\asciiabstract#1{\long\def\theasciiabstract{#1}}
\def\asciikeywords#1{\def\theasciikeywords{#1}}

\let\\\par\let\thelognumber\relax\let\thevolumenumber\relax
\let\thepapernumber\relax\let\thevolumeyear\relax\let\startpage\relax
\let\finishpage\relax\let\publishdate\relax\let\receiveddate\relax
\let\reviseddate\relax\let\accepteddate\relax\let\theasciititle\relax
\let\theasciiauthors\relax\let\theasciiaddress\relax
\let\theasciiabstract\relax\let\theasciikeywords\relax

\let\theasciiemail\relax

\ifplaintex
\font\logobig=cmssbx10 scaled 3836
\font\logomed=cmssbx10 scaled 2557
\else
\font\logobig=cmssbx10 scaled 4200
\font\logomed=cmssbx10 scaled 2800
\fi

\long\def\makeagttitle{   
\count0=\startpage
\agt\hfill      
\hbox to 45truept{\vbox to 0pt{\vglue -13truept{\logomed A\kern -.37em{\logobig 
T}\kern -.38em G}\vss}\hss}
\break
{\small Volume \thevolumenumber\ (\thevolumeyear)
\startpage--\finishpage\nl
Published: \publishdate}

\vglue .25truein

{\parskip=0pt\leftskip 0pt plus
1fil\def\\{\par\smallskip}{\Large\bf\thetitle}\par\medskip} \vglue
0.05truein

{\parskip=0pt\leftskip 0pt plus 1fil\def\\{\par}{\sc\theauthors}
\par\medskip}%
 
\vglue 0.03truein 

{\small\leftskip 25truept\rightskip 25truept{\bf Abstract}\stdspace\theabstract

{\bf AMS Classification}\stdspace\theprimaryclass
\ifx\thesecondaryclass\relax\else; \thesecondaryclass\fi\par
{\bf Keywords}\stdspace \thekeywords\par}\vglue 7truept

}   

\ifplaintex
\font\phead=cmsl9 scaled 950
\font\pnum=cmbx10 scaled 913
\font\pfoot=cmsl9 scaled 950
\headline{\vbox to 0pt{\vskip -4.5mm\line{\small\phead\ifnum
\count0=\startpage ISSN 1472-2739 (on-line) 1472-2747 (printed)
\hfill {\pnum\folio}\else\ifodd\count0\def\\{ }%
\ifx\theshorttitle\relax\thetitle\else\theshorttitle\fi\hfill{\pnum\folio}
\else\def\\{ and }{\pnum\folio}\hfill\ifx\theshortauthors\relax\theauthors
\else\theshortauthors\fi\fi\fi}\vss}}
\footline{\vbox to 0pt{\vglue 0mm\line{\small\pfoot\ifnum\count0=\startpage
\copyright\ \gtp\hfill\else
\agt, Volume \thevolumenumber\ (\thevolumeyear)\hfill\fi}\vss}}
\else
\font=cmsl9 scaled 1050
\font\lnum=cmbx10 
\font=cmsl9 scaled 1050
\makeatletter
\def\@oddhead{{\small\lhead\ifnum\count0=\startpage ISSN 1472-2739 
(on-line) 1472-2747 (printed)\hfill {\lnum\number\count0}\else\ifodd\count0
\def\\{ }\ifx\theshorttitle\relax \thetitle \else\theshorttitle\fi\hfill
{\lnum\number\count0}\else\def\\{ and }{\lnum\number\count0}
\hfill\ifx\theshortauthors\relax 
\theauthors\else\theshortauthors\fi\fi\fi}}\def\@evenhead{\@oddhead}
\def\@oddfoot{\small\lfoot\ifnum\count0=\startpage\copyright\ \gtp\hfill\else
\agt, Volume \thevolumenumber\ (\thevolumeyear)\hfill\fi}
\def\@evenfoot{\@oddfoot}
\makeatother
\fi
\let\maketitlepage\makeagttitle
\let\makeshorttitle\maketitlepage
\let\maketitle\maketitlepage

\newwrite\gtoutfile
\long\gdef\makeheadfile{  
{\def\\{, }\def\s{ }
\immediate\openout\gtoutfile head.xxx
\immediate\write\gtoutfile{To: math@arxiv.org}
\immediate\write\gtoutfile{Subject: put OR rep NNNNN:ppppp}
\immediate\write\gtoutfile{--text follows this line--}
\immediate\write\gtoutfile{Proxy-for: \ifx\theasciiauthors\relax
\theauthors\else\theasciiauthors\fi\s<\ifx\theasciiemail\relax\theemail\else\theasciiemail\fi>}
\immediate\write\gtoutfile{\noexpand\\}
\immediate\write\gtoutfile{Authors: \ifx\theasciiauthors\relax
\theauthors\else\theasciiauthors\fi}
{\def\\{ }\immediate\write\gtoutfile{Title: \ifx\theasciititle\relax
\thetitle\else\theasciititle\fi}}
\immediate\write\gtoutfile{Subj-class: GT or SG, GR etc}
\immediate\write\gtoutfile{MSC-class: \theprimaryclass\ifx\thesecondaryclass\relax\else, \thesecondaryclass\fi}
\immediate\write\gtoutfile{Journal-ref: Algebr. Geom. Topol. \thevolumenumber\s
(\thevolumeyear) \startpage-\finishpage}
\immediate\write\gtoutfile{Comments: Published by Algebraic and
Geometric Topology at}
\immediate\write\gtoutfile{\s\s\s  http://www.maths.warwick.ac.uk/agt/AGTVol\thevolumenumber/agt-\thevolumenumber-\thepapernumber.abs.html}
\immediate\write\gtoutfile{\noexpand\\}
\immediate\write\gtoutfile{}
\ifx\theasciiabstract\relax
\immediate\write\gtoutfile{\theabstract}\else
\immediate\write\gtoutfile{\theasciiabstract}\fi
\immediate\write\gtoutfile{}
\immediate\write\gtoutfile{\noexpand\\}
\immediate\write\gtoutfile{}
\immediate\closeout\gtoutfile}}  

\def\maketitlepage{\makeagttitle\makeheadfile}
\let\makeshorttitle\maketitlepage
\let\maketitle\maketitlepage

\def\ifplaintex{\expandafter\ifx\csname documentclass\endcsname\relax}

\def\gtp{{\mathsurround=0pt\it $\cal G\mskip-2mu$eometry \&\ 
$\cal T\!\!$opology $\cal P\!$ublications}}  

\def\recd{{\small Received:\qua\receiveddate\ifx\reviseddate\relax
\else\qquad Revised:\qua\reviseddate\fi\par}} 

\def\lognumber#1{\def\thelognumber{#1}}
\def\volumenumber#1{\def\thevolumenumber{#1}}
\def\volumeyear#1{\def\thevolumeyear{#1}}
\def\papernumber#1{\def\thepapernumber{#1}}
\def\pagenumbers#1#2{\def\startpage{#1}\def\finishpage{#2}}
\def\published#1{\def\publishdate{#1}}

\def\received#1{\def\receiveddate{#1}}
\def\revised#1{\def\reviseddate{#1}}
\def\accepted#1{\def\accepteddate{#1}}
\def\asciititle#1{\def\theasciititle{#1}}

\def\asciiaddress#1{\def\theasciiaddress{#1}}

\long\def\asciiabstract#1{\long\def\theasciiabstract{#1}}
\def\asciikeywords#1{\def\theasciikeywords{#1}}

\let\\\par\let\thelognumber\relax\let\thevolumenumber\relax
\let\thepapernumber\relax\let\thevolumeyear\relax\let\startpage\relax
\let\finishpage\relax\let\publishdate\relax\let\receiveddate\relax
\let\reviseddate\relax\let\accepteddate\relax\let\theasciititle\relax
\let\theasciiauthors\relax\let\theasciiaddress\relax
\let\theasciiabstract\relax\let\theasciikeywords\relax

\let\theasciiemail\relax

\ifplaintex
\font\logobig=cmssbx10 scaled 3836
\font\logomed=cmssbx10 scaled 2557
\else
\font\logobig=cmssbx10 scaled 4200
\font\logomed=cmssbx10 scaled 2800
\fi

\long\def\makeagttitle{   
\count0=\startpage
\agt\hfill      
\hbox to 45truept{\vbox to 0pt{\vglue -13truept{\logomed A\kern -.37em{\logobig 
T}\kern -.38em G}\vss}\hss}
\break
{\small Volume \thevolumenumber\ (\thevolumeyear)
\startpage--\finishpage\nl
Published: \publishdate}

\vglue .25truein

{\parskip=0pt\leftskip 0pt plus
1fil\def\\{\par\smallskip}{\Large\bf\thetitle}\par\medskip} \vglue
0.05truein

{\parskip=0pt\leftskip 0pt plus 1fil\def\\{\par}{\sc\theauthors}
\par\medskip}%
 
\vglue 0.03truein 

{\small\leftskip 25truept\rightskip 25truept{\bf Abstract}\stdspace\theabstract

{\bf AMS Classification}\stdspace\theprimaryclass
\ifx\thesecondaryclass\relax\else; \thesecondaryclass\fi\par
{\bf Keywords}\stdspace \thekeywords\par}\vglue 7truept

}   

\ifplaintex
\font\phead=cmsl9 scaled 950
\font\pnum=cmbx10 scaled 913
\font\pfoot=cmsl9 scaled 950
\headline{\vbox to 0pt{\vskip -4.5mm\line{\small\phead\ifnum
\count0=\startpage ISSN 1472-2739 (on-line) 1472-2747 (printed)
\hfill {\pnum\folio}\else\ifodd\count0\def\\{ }%
\ifx\theshorttitle\relax\thetitle\else\theshorttitle\fi\hfill{\pnum\folio}
\else\def\\{ and }{\pnum\folio}\hfill\ifx\theshortauthors\relax\theauthors
\else\theshortauthors\fi\fi\fi}\vss}}
\footline{\vbox to 0pt{\vglue 0mm\line{\small\pfoot\ifnum\count0=\startpage
\copyright\ \gtp\hfill\else
\agt, Volume \thevolumenumber\ (\thevolumeyear)\hfill\fi}\vss}}
\else
\font=cmsl9 scaled 1050
\font\lnum=cmbx10 
\font=cmsl9 scaled 1050
\makeatletter
\def\@oddhead{{\small\lhead\ifnum\count0=\startpage ISSN 1472-2739 
(on-line) 1472-2747 (printed)\hfill {\lnum\number\count0}\else\ifodd\count0
\def\\{ }\ifx\theshorttitle\relax \thetitle \else\theshorttitle\fi\hfill
{\lnum\number\count0}\else\def\\{ and }{\lnum\number\count0}
\hfill\ifx\theshortauthors\relax 
\theauthors\else\theshortauthors\fi\fi\fi}}\def\@evenhead{\@oddhead}
\def\@oddfoot{\small\lfoot\ifnum\count0=\startpage\copyright\ \gtp\hfill\else
\agt, Volume \thevolumenumber\ (\thevolumeyear)\hfill\fi}
\def\@evenfoot{\@oddfoot}
\makeatother
\fi
\let\maketitlepage\makeagttitle
\let\makeshorttitle\maketitlepage
\let\maketitle\maketitlepage

\newwrite\gtoutfile
\long\gdef\makeheadfile{  
{\def\\{, }\def\s{ }
\immediate\openout\gtoutfile head.xxx
\immediate\write\gtoutfile{To: math@arxiv.org}
\immediate\write\gtoutfile{Subject: put OR rep NNNNN:ppppp}
\immediate\write\gtoutfile{--text follows this line--}
\immediate\write\gtoutfile{Proxy-for: \ifx\theasciiauthors\relax
\theauthors\else\theasciiauthors\fi\s<\ifx\theasciiemail\relax\theemail\else\theasciiemail\fi>}
\immediate\write\gtoutfile{\noexpand\\}
\immediate\write\gtoutfile{Authors: \ifx\theasciiauthors\relax
\theauthors\else\theasciiauthors\fi}
{\def\\{ }\immediate\write\gtoutfile{Title: \ifx\theasciititle\relax
\thetitle\else\theasciititle\fi}}
\immediate\write\gtoutfile{Subj-class: GT or SG, GR etc}
\immediate\write\gtoutfile{MSC-class: \theprimaryclass\ifx\thesecondaryclass\relax\else, \thesecondaryclass\fi}
\immediate\write\gtoutfile{Journal-ref: Algebr. Geom. Topol. \thevolumenumber\s
(\thevolumeyear) \startpage-\finishpage}
\immediate\write\gtoutfile{Comments: Published by Algebraic and
Geometric Topology at}
\immediate\write\gtoutfile{\s\s\s  http://www.maths.warwick.ac.uk/agt/AGTVol\thevolumenumber/agt-\thevolumenumber-\thepapernumber.abs.html}
\immediate\write\gtoutfile{\noexpand\\}
\immediate\write\gtoutfile{}
\ifx\theasciiabstract\relax
\immediate\write\gtoutfile{\theabstract}\else
\immediate\write\gtoutfile{\theasciiabstract}\fi
\immediate\write\gtoutfile{}
\immediate\write\gtoutfile{\noexpand\\}
\immediate\write\gtoutfile{}
\immediate\closeout\gtoutfile}}  

\def\maketitlepage{\makeagttitle\makeheadfile}
\let\makeshorttitle\maketitlepage
\let\maketitle\maketitlepage

\lognumber{10}
\volumenumber{1}
\volumeyear{2001}
\papernumber{10}
\pagenumbers{201}{230}
\received{17 October 2000}
\revised{16 March 2001}
\accepted{13 April 2001}
\published{14 April 2001}

\usepackage{amscd}
\usepackage{amssymb}
\usepackage{amsmath}

\newtheorem{thm}{Theorem}[section]
\newtheorem{lem}[thm]{Lemma}
\newtheorem{theo}{Theorem}[section]
\newtheorem{coro}{Corollary}[section]
\newtheorem{remark}[thm]{Remark}
\newtheorem{cor}[thm]{Corollary}
\newtheorem{prop}[thm]{Proposition}

\theoremstyle{definition}
\newtheorem{example}[thm]{Example}
\newtheorem{defn}[thm]{Definition}

\newcommand{\GS}{\Gamma \sset}

\newcommand{\pr}{\operatorname{pr}}
\newcommand{\U}{{\mathcal{U}}}

\newcommand{\wedgel}{{\wedge^L}}
\newcommand{\op}{\operatorname{op}}
\newcommand{\R}{{\mathbb{R}}}
\newcommand{\Z}{{\mathbb{Z}}}

\newcommand{\C}{{\mathcal{C}}}
\newcommand{\N}{{\mathbb{N}}}

\newcommand{\gr}{\operatorname{gr}}

\newcommand{\sd}{\operatorname{sd}}
\newcommand{\trf}{\operatorname{trf}}
\newcommand{\HH}{\operatorname{HH}}
\newcommand{\THH}{\operatorname{THH}}

\newcommand{\pMap}{\operatorname{Map}_*}
\newcommand{\Tor}{\operatorname{Tor}}

\newcommand{\id}{\operatorname{id}}

\newcommand{\sset}{{\mathcal{S}_*}}

\newcommand{\holim}[1]{\mathop{\underset{#1}
            {{\text{\rm holim}}}}}
\newcommand{\colim}[1]{\mathop{\underset{#1}
            {{\text{\rm colim}}}}} 
\newcommand{\hocolim}[1]{\mathop{\underset{#1}
            {{\text{\rm hocolim}}}}}

\begin{document}

\title[Filtered $TC$ and $K$-theory of nilpotent ideals]{Filtered
  Topological Cyclic Homology \\ and relative 
  $K$-theory of nilpotent ideals}
\asciititle{Filtered Topological Cyclic Homology and relative K-theory
of nilpotent ideals}
\author{Morten Brun}
\email{brun@math.u-strasbg.fr}
\address{Institut de Recherche
Mathematique Avanc\'ee\\
CNRS et Universit\'e Louis Pasteur, 7 rue R. Descartes\\
67084 Strasbourg Cedex, France}
\asciiaddress{Institut de Recherche
Mathematique Avancee\\
CNRS et Universite Louis Pasteur, 7 rue R. Descartes\\
67084 Strasbourg Cedex, France}

\begin{abstract}
  In this paper certain filtrations of topological Hochschild homology
  and topological cyclic homology are examined. As an example we show
  how the filtration with respect to a nilpotent ideal gives rise to
  an analog of a theorem of Goodwillie saying that rationally relative
  $K$-theory and relative cyclic homology agree. Our variation says
  that the $p$-torsion parts agree in a range of degrees. We use it to
  compute $K_i(\Z/p^n)$ for $i \le p-3$.
\end{abstract}

\asciiabstract{
  In this paper we examine certain filtrations of topological
  Hochschild homology and topological cyclic homology. As an example
  we show how  
  the filtration with respect to a nilpotent ideal gives rise to
  an analog of a theorem of Goodwillie saying that rationally
  relative K-theory and relative cyclic homology agree. Our
  variation says that the
  p-torsion parts agree in a range of degrees. We use it to compute
  K_i(Z/p^m) for i < p-2.}

\asciikeywords{K-theory, topological Hochschild homology, cyclic
homology, topological cyclic homology}

\maketitle

\section{Introduction}

The aim of this paper is to examine a certain filtration of
topological Hochschild homology of a functor with smash product
equipped with a filtration. The former filtration preserves the cyclic
structure and it induces a filtration of topological cyclic homology.
By a theorem of McCarthy \cite{McCarthy} topological cyclic homology
is closely related to algebraic $K$-theory, and in some interesting
cases topological cyclic homology determines the $K$-groups.  The
methods developed in this paper stem from a paper of Hesselholt and
Madsen, where the $K$-groups for finite algebras over Witt vectors of
perfect fields of positive characteristic are computed
\cite{Hesselholt-Madsen}. One difference is that here general
filtrations are considered, while the filtrations considered in
\cite{Hesselholt-Madsen} are split.  In the paper \cite{Brun} the
filtration of $\Z/p^n$ by the powers of the ideal $p\Z/p^n$ was used
to compute topological Hochschild homology of the ring $\Z/p^n$.  This
is the example that motivated the generality of the present paper.

Given a ring $R$
with an ideal $I$, we shall let $K(R,I)$ denote the homotopy fibre of
the map 
$K(R) \rightarrow K(R/I)$, and we shall let $HC(R,I)$ denote the
homotopy fibre
of the map $HC(R) \rightarrow HC(R/I)$.
As an example of how the filtrations constructed can be useful, we
prove the following analog of a theorem of Goodwillie
\cite{Goodwillie}.

\begin{theo}
  Let $R$ be a simplicial ring with an ideal $I$ satisfying $I^m =
  0$. Suppose that $R$ and $R/I$ are flat as modules over $\Z$. Then
  there is an isomorphism of homotopy
  groups of $p$-adic completions
  \begin{displaymath}
    \pi_i K(R,I)^{\wedge}_p \cong \pi_{i-1} HC(R,I)^{\wedge}_p 
  \end{displaymath}
  when $0 \le i < p/(m-1) -2$ and a surjection
  \begin{displaymath}
    \pi_i K(R,I)^{\wedge}_p \rightarrow \pi_{i-1} HC(R,I)^{\wedge}_p 
  \end{displaymath}
  when $i < p/(m-1) -1$.  
\end{theo}

In the case where $R$ and $R/I$ are not flat as modules over $\Z$, we
can replace them with weakly equivalent simplicial rings that are
degreewise free abelian groups. Since $K(R,I)$ is homotopy invariant
we obtain that $K_i(R,I)^{\wedge}_p$ is isomorphic to the $p$-adic
completion of derived relative cyclic homology in the same range of degrees.
In section \ref{derivedcyclic} we recall the definition of derived
cyclic homology, 
and we 
compute enough derived cyclic homology groups for $\Z/p^n$ to deduce
the following result. 
\begin{coro}
  For $1 \le i \le p-3$, the $K$-groups of $\Z/p^n$ are:
  \begin{displaymath}
     \pi_i K(\Z/p^n) \cong
    \begin{cases}
      0 & \text{if $i$ is even} \\
      \Z/p^{j(n-1)} (p^j -1) & \text{if $i =2j-1$}
    \end{cases}
  \end{displaymath}
\end{coro}
The starting point of the above result is Quillen's computation of
$\pi_*K(\Z/p)$ in \cite{QuillenK}.
The result agrees with the computation of $K_i(\Z/p^n)$ for $0 \le i \le 4$
of Aisbett, Puebla and Snaith \cite{aisbett-puebla-snaith} starting from 
Evens and Friedlander's computation 
of $K_i(\Z/p^2)$ for $0 \le i \le 4$ and for $p \ge 5$
\cite{Evens-Friedlander}. It also shows that the homotopy groups of
$BGL(\Z/p^n)^+$ and of the homotopy fibre of $\psi^{p^n} - \psi^{p^{n-1}} : BU
\rightarrow BU$ are different so these spaces can not be homotopy
equivalent, as was also proven by Priddy in
\cite{Priddy-on-conjecture} in the case $n=2$.

In view of Quillen's computation of $\pi_* K(\Z/p)$ only the
$p$-torsion part of corollary \ref{kofzpsquare} is hard to 
prove. Let us show that if $l$ is relatively prime to $p$ then the natural map
$K(\Z/p^n) \rightarrow K(\Z/p)$ induces an isomorphism on
homotopy groups with coefficients in $\Z/l$. Since $BGL(\Z/p^n)^+$ and
$BGL(\Z/p)^+$ are simple spaces it suffices by the Whitehead
theorem to show that the map $BGL(\Z/p^n)^+ \rightarrow BGL(\Z/p^{n-1})^+$
induces an isomorphism on homology with coefficients in $\Z/l$. 
For this we note that the kernel
of the map $GL_m(\Z/p^n) 
\rightarrow GL_m(\Z/p^{n-1})$ consists of matrices of the form
$I+p^{n-1}M$. Multiplication is given by $(I+p^{n-1}M)(I+p^{n-1}N) = I
+ p^{n-1}(M+N)$ so the kernel $J$ of the map $GL(\Z/p^n) \rightarrow
GL(\Z/p^{n-1})$ is a vectorspace over $\Z/p$. The Serre spectral
sequence 
\begin{displaymath}
  H_*(BGL(\Z/p^{n-1}),H_*(BJ,\Z/l)) \Rightarrow H_*(BGL(\Z/p^n),\Z/l)
\end{displaymath}
associated to the fibration $BJ \rightarrow BGL(\Z/p^n) \rightarrow
BGL(\Z/p^{n-1})$ collapses to an isomorphism $H_*(BGL(\Z/p^{n-1},\Z/l)
\cong H_*(BGL(\Z/p^n),\Z/l)$.

Only elementary properties of the filtrations of topological
Hochschild homology and topological cyclic homology are studied in this
note. The focus is on a filtered version of the norm cofibration
sequence for the fixed points of topological Hochschild homology.
Traditionally, for example in \cite{Bokstedt-Madsen} and in
\cite{Hesselholt-Madsen}, the role of the norm cofibration sequence is that it
allows one to
determine the fixed point spectra inductively. Here we use it to keep
track of the connectivity properties of our filtration of topological
cyclic homology.

The paper is organized as follows: In section
\ref{filtered-Hochschild} generalities on filtrations of 
monoids in a symmetric monoidal category are given. It is noted that a
filtered monoid is a monoid in the symmetric monoidal category of
filtered objects, and 
therefore it fits into the Hochschild construction. In section
\ref{THH} a filtered functor with smash product is defined to be
a filtered monoid in the category of Gamma spaces, and fundamental
properties of the topological Hochschild homology of a filtered
functor with smash product are established. In section
\ref{cyclotomic} we introduce the concept of a cyclotomic filtered
Gamma space. This is a filtered Gamma space with an action of the circle
group having enough extra properties to make it possible to construct
a filtered version 
of topological cyclic homology out of it. It is shown that topological
Hochschild homology of 
a filtered functor with smash product is such a a cyclotomic filtered
Gamma space.
In section \ref{goodwillietype} a proof of
theorem \ref{mainthm} 
is given. In section \ref{derivedcyclic} we compute enough
derived cyclic homology of the ring $\Z/p^n$ to prove corollary
\ref{kofzpsquare}. 

It might be appropriate add a remark on terminology. Following
Bousfield and Friedlander \cite{Bousfield-Friedlander} we do not assume
Gamma spaces to be special, and we do not assume spectra to
to be omega-spectra. 

\rk{Acknowledgments} 
This work was funded by a Marie Curie
Fellowship. The author wishes to thank F. Waldhausen and
the University of Bielefeld for hospitality while this paper was
written and he wants to thank M. B\"okstedt for his
constant interest in and collaboration on the computations motivating
the present work. In particular B\"okstedt proved proposition
\ref{folklore} 
and proposition \ref{derivedcycliccom} before the author. His proofs
are different from the ones presented here. Finally the author wants to thank
the referee for a lot of helpful comments.

\section{The Filtered Hochschild Construction}
In this section we shall study filtered objects in a category
$\C$. 
\label{filtered-Hochschild}
\subsection{Filtered Objects}\label{definefiltered}
A {\em filtered object} in a category $\C$ is a functor from the category
$\Z$,  with exactly one morphism $n \rightarrow m$ if $n \le m$, to
$\C$. That is, a filtered object is a sequence
\begin{displaymath}
  \dots \rightarrow
  X({-1}) \rightarrow X(0) \rightarrow X(1) \rightarrow \dots \rightarrow
  X(n) \rightarrow X({n+1}) \rightarrow \dots
\end{displaymath}
of composable morphisms in $\C$.
A morphism of filtered objects is simply a natural
transformation. For some choices of $\C$ there is a functor $H$ from
the category $\C^\Z$ of filtered objects in $\C$ to the category of
exact couples of 
(graded) abelian groups in the sense of Massey \cite{Massey1952}.

\begin{example} Functors from filtered objects to exact couples:
  \begin{enumerate}
  \item The category of chain complexes and injective
    chain homomorphisms together with the functor $H$ given by
    homology.
  \item We can take $\C$ to be the category of topological spaces and
    cofibrations 
    and let $H$ be given by (generalized) homology.  
  \end{enumerate}
\end{example}
Given
objects $X_i$ in $\C$ we shall denote their coproduct by $\bigvee_i
X_i$, and given
a diagram $Z \leftarrow X \rightarrow Y$ we shall 
denote the colimit, that is, the pushout by $Z \cup_X Y$. 
\begin{lem}
\label{pushoutlemma}
  Given a functor $F : \Z \times \Z \rightarrow \C$ and $k \in \Z$ the
  following diagram is a pushout diagram:
  $$\begin{matrix} 
     \bigvee_{i+j=k} F(i-1,j) \cup_{F(i-1,j-1)} F(i,j-1) 
  &\to& 
     \colim{i+j \le  k-1} F(i,j) \\
  \downarrow && \downarrow \\
     \bigvee_{i+j=k} F(i,j) 
      & \to &  
  \colim{i+j\le k} F(i,j). \end{matrix}$$
\end{lem}

\subsection{Filtered objects in monoidal categories}\label{filteredmonoids}
From now on $\C = (\C,\otimes,I)$ shall denote a cocomplete
symmetric monoidal category.
Given filtered objects $X$ and $Y$ in $\C$ we can define a filtered
object $X \otimes Y$ in $\C$ by letting $(X \otimes Y)(k) 
= \colim{i+j \le k} X(i) \otimes Y(j)$. The pairing $\otimes : \C^\Z
\times \C^\Z \rightarrow \C^\Z$ defines a symmetrical monoidal
structure on $\C^\Z$ with
unit $I$ given by the filtered object with $I(k)$ equal to
the initial object in $\C$ for $k<0$ and with $I(k)$ equal to the unit for
the monoidal structure of $\C$ when $k \ge 0$. We have consciously
chosen the same symbols for the pairing and unit in $\C^\Z$ as in $\C$
because we can consider $\C$ as a full symmetric monoidal subcategory of
$\C^\Z$. 

A {\em filtered monoid} in $\C$ is a monoid in the category
$\C^\Z$. Explicitly, a filtered monoid in $\C$ is a sequence
\begin{displaymath}
  \dots \rightarrow M(-1) \rightarrow M(0) \rightarrow M(1)
  \rightarrow \dots \rightarrow M(n) \rightarrow M(n+1) \rightarrow \cdots
\end{displaymath}
of composable morphisms in $\C$ together with morphisms
\begin{gather*}
  \mu_{i,j} : M(i) \otimes M(j) \rightarrow M(i+j) \\
  \eta : I \rightarrow M(0),
\end{gather*}
satisfying the following relations for associativity and unitality:
\begin{gather*}
  \mu_{i+j,k} \circ (\mu_{i,j} \otimes \id_{M(k)}) = \mu_{i,j+k} \circ
  (\id_{M(i)} \otimes \mu_{j,k}), \\ 
  \mu_{0,i} \circ (\eta \otimes \id_{M(i)}) =
  \lambda_{M(i)}, \\ 
  \mu_{i,0} \circ (\id_{M(i)} \otimes \eta) = \rho_{M(i)}. 
\end{gather*} 
Here $\lambda_{M(i)} : I \otimes M(i) \cong M(i)$ and $\rho_{M(i)} : M(i)
\otimes I \cong M(i)$ are part of the symmetric monoidal structure of $\C$.
We shall call a filtered monoid in the category of abelian groups a {\em
  filtered ring}.

If $*$ is a terminal object of $\C$ and $X \rightarrow Y$ is a map in
$\C$, we shall denote any choice of pushout of the diagram $*
\leftarrow X \rightarrow Y$ by $Y/X$. We shall say that the product
$\otimes$ of $\C$ commutes with quotients if there is a natural isomorphism
$({X_1}/{X_2}) \otimes Y \cong ({X_1 \otimes Y})/({X_2 \otimes Y})$.
\begin{lem}\label{idquotients}
  If $\C$ is a cocomplete symmetric monoidal category with a terminal
  object, then 
  given filtered objects $X$ and $Y$ of $\C$
  there is an isomorphism:
  \begin{displaymath}
    \frac{(X\otimes Y)(k)}{(X \otimes Y)(k-1)} \cong \bigvee_{i+j
    =k} \frac{X(i) \otimes Y(j)}{X(i-1)\otimes Y(j) \cup_{X(i-1)
    \otimes Y(j-1)} X(i) \otimes Y(j-1)}.
  \end{displaymath}
  If in addition $\otimes$ commutes with quotients, then there is an
  isomorphism: 
  \begin{displaymath}
    \frac{(X\otimes Y)(k)}{(X \otimes Y)(k-1)} \cong \bigvee_{i+j
    =k} \frac{X(i)}{X(i-1)} \otimes \frac{Y(j)}{Y(j-1)}.
  \end{displaymath}
\end{lem}
\proof
  For the first part, it suffices to note that by lemma
  \ref{pushoutlemma} the following diagram
  in $\C$ is a pushout: 
  $$\begin{matrix} 
        \bigvee_{i+j
      =k} {X(i-1)\otimes Y(j) \cup_{X(i-1)
      \otimes Y(i-1)} X(i) \otimes Y(j-1)}
  &\to& 
  (X \otimes Y)(k-1) \\
  \downarrow && \downarrow \\
      \bigvee_{i+j=k} X(i) \otimes Y(j) & \to & 
  (X \otimes Y)(k). \end{matrix}$$
  For the second part, we note that:
  \begin{displaymath}
    ({X_1}/{X_0}) \otimes ({Y_1}/{Y_0}) \cong
    \frac{({X_1}/{X_0}) \otimes Y_1}{({X_1}/{X_0}) \otimes Y_0}
    \cong \frac{({X_1 \otimes Y_1})/({X_0 \otimes Y_1})}{({X_1
    \otimes Y_0})/({X_0 \otimes Y_0})}
  \end{displaymath}
  for $X_0 \rightarrow X_1$ and $Y_0 \rightarrow Y_1$ maps in $\C$,
  and that
  given a map $B \cup_D C \rightarrow A$ in $\C$ we have:
  $$\frac{A/B}{C/D} \cong \frac{A}{B \cup_{D} C}.\eqno{\qed}$$

\subsection{The Hochschild construction}
Let $\C$ denote a symmetric monoidal category. The
Hochschild construction is a functor $Z$ from the category of monoids
in $\C$ to Connes' category of cyclic objects in $\C$. A good
reference for the category of cyclic objects is the book of Loday
\cite[chapter 6]{Loday}. Given a monoid
$M$ in $\C$, $Z(M)$ is defined as follows: It has $n$-simplices
\begin{displaymath}
  Z_n(M) = M \otimes \dots \otimes M \qquad \text{$(n+1)$ factors.}
\end{displaymath}
The cyclic operator is given by the automorphism $t_n$ of 
$M \otimes \dots \otimes M$ 
cyclically shifting the $(n+1)$ factors to the right.
The face maps are given by the formula:
\begin{displaymath}
  d_i = t_{n-1}^{i} \circ (\mu \otimes \id) \circ t_{n}^{-i}, \qquad{0
  \le i \le n},
\end{displaymath}
where $\mu : M\otimes M \rightarrow M$ is the multiplication in $M$.
The degeneracies
are given by the formula: 
\begin{displaymath}
  s_i = t_{n+1}^{(i+1)} \circ (\eta \otimes \id) \circ t^{-({i+1})}_{n},
  \qquad 0 \le i \le n,
\end{displaymath}
where $\eta : I \rightarrow M$ is the unit in $M$.

Since $\C^\Z$ is equipped with a symmetric monoidal structure, we can also consider the
Hochschild construction on 
monoids in $\C^\Z$, that is on filtered monoids in $\C$.
\begin{prop}\label{Zquotients}
  Let $M$ be a filtered monoid in a cocomplete symmetric monoidal
  category $\C$, where $\otimes$ commutes with quotients. Then for
  each $n \ge 0$ there is an isomorphism of cyclic objects
  \begin{displaymath}
    \frac{Z_n(M)(k)}{Z_n(M)(k-1)} \cong \bigvee_{i_0 + \dots + i_n =
    k} \frac{M(i_0)}{M(i_0-1)} \otimes \dots \otimes \frac{M(i_n)}{M(i_n-1)}.
  \end{displaymath}
\end{prop}
\begin{proof}
  This is a direct consequence of lemma \ref{idquotients}.
\end{proof}
The above proposition can be reformulated in terms of the associated
graded monoid $\gr M$ for $M$. Here $\gr M$ is the filtered monoid in $\C$
with
\begin{displaymath}
  \gr(M) (k) = \bigvee_{i \le k} M(i)/M(i-1),
\end{displaymath}
and with multiplication induced by the maps
\begin{align*}
  \frac{M(i)}{M(i-1)} \otimes \frac{M(j)}{M(j-1)} &\cong \frac{M(i)
  \otimes M(j)}{M(i) \otimes M(j-1) \cup_{M(i-1) \otimes M(j-1)}
  M(i-1) \otimes M(j)} \\
  &\rightarrow \frac{M(i+j)}{M(i+j-1)}.
\end{align*}
The proposition says that the filtration quotients for $Z(M)$
and $Z(\gr M)$ are isomorphic.

\section{Filtered Topological Hochschild Homology}
\label{THH}
\subsection{Topological Hochschild homology}
\label{defTHH}
We briefly recall the definition of  topological Hochschild homology:
Let $I$ denote the category with one object $n$ for each integer $n
\ge 0$ and with $I(m,n)$ given by the set of injective maps from $\{1,
\dots m \}$ to $\{1,\dots ,n\}$. 
Let $L$ denote 
a functor with smash product in the sense of B\"okstedt
\cite{Bokstedt1} or in the more restrictive sense described
below. $\THH (L)$ is the cyclic pointed 
simplicial set with
$k$-simplices equal to the homotopy colimit
$$\hocolim{(i_0,\dots,i_k) \in I^{k+1}} F (S^{i_0} \wedge \dots \wedge
S^{i_k}, L(S^{i_0}) 
\wedge \dots \wedge L(S^{i_k}))$$
and with structure maps of the same type as for the
Hochschild construction.
Details on this construction can be found in \cite{BHM}.
The symbol $F$ denotes derived
function space, that is, if $X$ and $Y$ are pointed simplicial sets,
then $  F(X,Y) = \sset(X,\sin |Y|)$,
where $\sin|Y|$ denotes the singular complex on the geometric
realization of $Y$, and $\sset$ denotes the internal function object
in the category $\sset$ of pointed simplicial sets. Occasionally we
shall write $\Omega^n Y = F(S^n,Y)$ for the $n$'th loop space of $Y$.

For the purpose of this note, a functor with smash product is a
monoid in the category $\GS$ of Gamma spaces considered for example by
Bousfield
and Friedlander \cite[definition 3.1]{Bousfield-Friedlander}.  
Let us recall that a Gamma space is a pointed functor from the
category $\Gamma$ of pointed finite sets to the category $\sset$ of
pointed simplicial sets. To be precise
$\Gamma$ is
the category with one object $n^+ = \{ 0,1,
\dots, n\}$ for each $n\ge 0$, and with $\Gamma(m^+,n^+)$ the set functions
from $m^+$ to $n^+$ fixing $0$. 
A pointed category is a category with an object which is both initial
and terminal, and a functor between pointed categories is pointed if
it takes an object which is both initial and final to an object of the
same kind. 
Let us stress that our notion of a Gamma space is different from the
notion in the paper by Segal \cite{Segal}. 
Given two Gamma spaces $X$ and $Y$, their smash product is the
Gamma space $X \wedge Y$ with 
\begin{displaymath}
  (X \wedge Y)(n^+) = \colim{n_1^+ \wedge n_2^+ \rightarrow n^+}
  X(n_1^+ ) \wedge Y(n_2^+).
\end{displaymath}
The unit for the operation $\wedge$ is the functor $\mathbb S : \Gamma
\rightarrow \sset$ with
$\mathbb S(n^+) = n^+$. Lydakis noted in \cite[Theorem 2.18]{Lydakis} that the
category of Gamma spaces 
is a symmetric monoidal category with respect to the smash product
pairing and unit. 
By definition a functor with smash product $L$ is a monoid in the category
$\GS$. Explicitly this means that $L$ is a pointed functor $L :
\Gamma \rightarrow \sset$ together with natural transformations 
\begin{gather*}
  \mu: L(m^+) \wedge L(n^+) \rightarrow L(m^+ \wedge n^+), \\
  \eta: n^+ \rightarrow L(n^+),
\end{gather*}
satisfying the following relations for associativity and unitality:
\begin{gather*}
  \mu \circ (\mu \wedge \id) = \mu \circ (\id \wedge \mu), \quad
  \mu \circ (\eta \wedge \id) = \lambda, \quad
  \mu \circ (\id \wedge \eta) = \rho,
\end{gather*}
where $\lambda :m^+ \wedge L(n^+) \rightarrow L(m^+ \wedge n^+)$ is
adjoint to the map
\begin{displaymath}
  m^+ \rightarrow \Gamma (n^+, m^+ \wedge n^+) \overset{L}
  \rightarrow 
  \sset(L(n^+),L(m^+ \wedge n^+))
\end{displaymath}
and $\rho : L(m^+) \wedge n^+ \rightarrow L(m^+ \wedge n^+)$ is
adjoint to the map
\begin{displaymath}
  n^+ \rightarrow \Gamma (m^+, m^+ \wedge n^+) \overset{L}
  \rightarrow \sset(L(m^+),L(m^+ \wedge n^+)).
\end{displaymath}

\begin{example}
  For this note the most relevant example of an FSP is the functor
  $\widetilde \Z : \Gamma \rightarrow \sset$ with $\widetilde \Z (n^+)
  = \Z\{n^+ \}/\Z \{0\}$ the reduced free abelian group on the pointed
  set $n^+ = \{ 0,1,\dots,n\}$. The multiplication $\mu$ is given by the
  composition 
  \begin{displaymath}
    \mu: \widetilde \Z(m^+) \wedge \widetilde \Z(n^+) \rightarrow
  \widetilde \Z(m^+) \otimes_{\Z} \widetilde \Z(n^+) \cong \widetilde
  \Z(m^+ \wedge n^+),
  \end{displaymath}
  and the unit is given by the inclusion of the basis $n^+$ in the
  free abelian group $\Z \{n^+ \}$ composed with the quotient map.
  Given any ring $R$ we obtain an FSP $\widetilde R$ with $\widetilde R
  (n^+) = R \otimes \widetilde \Z(n^+)$. The multiplication and the
  unit in $\widetilde R$ are explained in example
  \ref{filtringfiltfsp} below.
\end{example}
Given a Gamma space $X$, we can extend it to a
functor $X_1$ defined on the category of pointed sets by letting 
\begin{displaymath}
  X_1(K) = \colim{n^+ \rightarrow K} X(n^+),
\end{displaymath}
for $K$ a pointed set, 
and we can define an endofunctor $X_2$ on $\sset$ by letting
$(X_2(U))_k = (X_1(U_k))_k$ for $U$ a pointed simplicial set. From now
on we shall not distinguish notationally between a Gamma space and the
induced endofunctor on $\sset$.

Given a Gamma space $X$ and pointed simplicial sets $U$ and $V$, there is
a map $X(U) \wedge V \rightarrow X(U \wedge V)$ obtained by applying
the above map 
$\rho$ degreewise.
The following lemma is given in \cite[prop. 5.21]{Lydakis}:
\begin{lem}\label{gammaconnected}
  If $U$ is $m$-connected and $V$ is $n$-connected, then the map $X(U)
  \wedge V \rightarrow X(U \wedge V)$ is $2m+n+3$-connected.
\end{lem}
Together with the approximation lemma of B\"okstedt 
(see either \cite{Bokstedt1} or \cite[lemma 2.5.1]{Brun}), it
can be used to
prove the following.
\begin{lem}\label{approxlemma}
  Given a Gamma space $X$ and $(j_0,\dots,j_k) \in I^{k+1}$, then the map
\begin{multline*}
F (S^{j_0} \wedge \dots \wedge
S^{j_k}, X(S^{j_0}) \wedge \dots \wedge X(S^{j_k})) \rightarrow \\
\hocolim{(i_0,\dots,i_k) \in I^{k+1}} F (S^{i_0} \wedge \dots \wedge
S^{i_k}, X(S^{i_0}) \wedge \dots \wedge X(S^{i_k}))
\end{multline*}
is $j-1$-connected. Here $j$ denotes the minimum of the cardinalities of 
$j_0, \dots, j_k$.
\end{lem}
Given an FSP $L$ and a finite pointed set $n^+$, we shall let
$\THH(L;n^+)$ denote the cyclic pointed simplicial set with
$k$-simplices equal to the homotopy colimit
\begin{displaymath}
  \hocolim{(i_0,\dots,i_k) \in I^{k+1}}F(S^{i_0} \wedge \dots \wedge S^{i_k}, L(S^{i_0})
  \wedge \dots \wedge L(S^{i_k}) \wedge n^+),
\end{displaymath}
where $n^+$ acts as a dummy variable for the cyclic structure.
There is an endofunctor $\THH(L,-)$ on $\sset$ associated to the
Gamma space $n^+ \mapsto \THH(L;n^+)$. We shall freely use the identification
$\THH(L) \cong \THH(L;1^+) \cong \THH(L;S^0)$. 
\begin{lem}\label{thhisomega}
  The map $\THH(L) \rightarrow \Omega^n \THH(L;S^n)$, adjoint to 
  $\THH(S;S^0) \wedge S^n \rightarrow \THH(L ; S^n)$, is a weak
  equivalence. 
\end{lem}
\begin{proof}
  By the work of Segal \cite[prop. 1.4]{Segal} it suffices to show that the
  Gamma space $n^+ \mapsto \THH(L;n^+)$ is {\em very special}, that is,
  the map 
  $$(\pr_{1*},\pr_{2*}):\THH(L;m^+ \vee n^+) \rightarrow \THH(L;m^+) \times
  \THH(L;n^+)$$ 
  induced by the projections $\pr_{1} : m^+ \vee n^+ \rightarrow m^+ 
  \vee 0^+ = m^+$ and $\pr_{2} : m^+ \vee n^+ \rightarrow 0^+ \vee n^+  
  = n^+$,
  is a weak equivalence, and that the monoid $\pi_0 \THH(L;1^+)$
  is a group. By lemma \ref{shearingtrick} below it suffices to show
  that the Gamma spaces $n^+ \mapsto \THH_k(L;n^+)$ are very
  special. To see that the map 
  $$\THH_k(L;m^+ \vee n^+) \rightarrow \THH_k(L;m^+) \times
  \THH_k(L;n^+)$$ 
  is a weak equivalence, it suffices by the approximation lemma
  \ref{approxlemma} to note that by the Whitehead theorem the map
  \begin{align*}
    &L(S^{i_0}) \wedge \dots \wedge L(S^{i_k}) \wedge (m^+ \vee n^+) \\
    & \cong (L(S^{i_0}) \wedge \dots \wedge L(S^{i_k}) \wedge m^+) \vee
    (L(S^{i_0}) \wedge \dots \wedge L(S^{i_k}) \wedge n^+) \\
    & \rightarrow  (L(S^{i_0}) \wedge \dots \wedge L(S^{i_k}) \wedge
    m^+) \times 
    (L(S^{i_0}) \wedge \dots \wedge L(S^{i_k}) \wedge n^+) \\
  \end{align*}
  is $2(i_0 + \dots + i_k)-1$-connected.

  To see that $\pi_0 \THH_k(L;1^+)$ is a group, it suffices to note
  that $\pi_0 F(S^{i_0} \wedge \dots \wedge S^{i_k}, L(S^{i_0}) \wedge
  \dots \wedge L(S^{i_k}))$ is a group.
\end{proof}
We owe the following lemma to S. Schwede.
\begin{lem}\label{shearingtrick}
  Let $X$ be a simplicial Gamma space, and assume that for each $k$,
  $X_k$ is a very special Gamma space. Then the Gamma topological
  space $|X|$ sending 
  $n^+$ to the realization of $[k] \mapsto X_k(n^+)$ is very
  special.
\end{lem}
\begin{proof}
  It follows from the realization lemma and the fact that realization
  commutes with products that the resulting Gamma space
  is special, that is, the map $|X(m^+ \vee n^+)| \rightarrow
  |X(m^+)| \times |X(n^+)|$ induced by the projections $\pr_{1}$ and
  $\pr_{2}$ is a homotopy equivalence. A special Gamma space $Y$ is very
  special
  when the monoid $\pi_0 |Y(1^+)|$ with 
  multiplication induced by the composite
  \begin{displaymath}
    |Y(1^+)| \times |Y(1^+)| \overset{f} \rightarrow |Y(2^+)|
    \overset{Y(\mu)} \rightarrow |Y(1^+)|
  \end{displaymath}
  is a group.
  Here $\mu: 2^+ \rightarrow 1^+$ is the fold map with $\mu(i) = 1$
  for $i = 1,2$ and $f$ is a homotopy inverse to the homotopy
  equivalence $|Y(2^+)| \rightarrow |Y(1^+)| \times
  |Y(1^+)|$.
  This is
  equivalent to the map $(Y(\mu), Y(\pr_2)): |Y(2^+)| \rightarrow
  |Y(1^+)| \times |Y(1^+)|$ being a homotopy equivalence. (Clearly if this
  map is a homotopy equivalence, then $\pi_0 Y(1^+)$ is a
  group. Conversely, if $\pi_0 Y(1^+)$ is a group, then this map
  induces an isomorphism on $\pi_n$ for all $n \ge 0$, and by the
  Whitehead theorem we can conclude that it is a homotopy equivalence.)
  It follows from
  the realization lemma that $|X|$ is very special.
\end{proof}
\subsection{Filtered Topological Hochschild Homology}
\label{filteredFSP} 
To make a filtered version of topological Hochschild homology we
replace the category $\sset$ of pointed simplicial sets by the
category $\sset^\Z$ of filtered 
pointed simplicial sets. By a {\em Gamma filtered space} we shall mean a 
pointed functor from $\Gamma$ to $\sset^\Z$. The smash product of two
Gamma filtered spaces $X$ and $Y$, given by the
formula 
\begin{displaymath}
  (X \wedge Y) (n^+) = \colim{n_1^+ \wedge n_2^+ \rightarrow n^+}
  X(n_1^+) \wedge Y(n_2^+),
\end{displaymath}
makes the category $\GS^\Z$ of Gamma filtered spaces into a symmetric
monoidal category.
A {\em filtered FSP} is a monoid in the category $\GS^\Z$. Explicitly a
filtered FSP can be described as a 
functor 
$L : \Gamma \times \Z \rightarrow \sset$ together with natural
transformations 
\begin{gather*}
  {\mu} : L(m^+,s) \wedge L(n^+,t)  \rightarrow L(m^+ \wedge
  n^+,s+t) \\
  {\eta} : n^+  \rightarrow L(n^+,0)
\end{gather*}
satisfying the following relations:
\begin{displaymath}
  \mu \circ (\mu \wedge \id) = \mu \circ (\id \wedge \mu), \quad \mu
  \circ (\eta \wedge \id) = \lambda, \quad \mu \circ (\id \wedge \eta)
  = \rho,
\end{displaymath}
where $\lambda :m^+ \wedge L(n^+,s) \rightarrow L(m^+ \wedge n^+,s)$ is
adjoint to the map
\begin{displaymath}
  m^+ \rightarrow \Gamma (n^+, m^+ \wedge n^+) \overset{L(-,s)}
  \rightarrow 
  \sset(L(n^+,s),L(m^+ \wedge n^+,s))
\end{displaymath}
and $\rho : L(m^+,s) \wedge n^+ \rightarrow L(m^+ \wedge n^+,s)$ is
adjoint to the map
\begin{displaymath}
  n^+ \rightarrow \Gamma (m^+, m^+ \wedge n^+) \overset{L(-,s)}
  \rightarrow \sset(L(m^+,s),L(m^+ \wedge n^+,s)).
\end{displaymath}

Note that the category of Gamma filtered spaces is isomorphic to the
category of filtered Gamma spaces, and hence
a filtered FSP also can be described as being a filtered
monoid in the category of Gamma spaces.
\begin{example}
  \label{filtringfiltfsp}
  Given a filtered ring $R$ (that is, a filtered monoid in the
  category of abelian groups) there is a filtered FSP $\widetilde R$
  with $\widetilde R (n^+,s) = \widetilde \Z(n^+) \otimes_{\Z} R(s)$. The
  multiplication is given by the composition
  \begin{align*}
    \mu : \widetilde R (m^+,s) \wedge \widetilde R (n^+,t)& \rightarrow
  \widetilde R (m^+,s) \otimes_{\Z} \widetilde R (n^+,t)
  \\ 
  &\cong 
  \widetilde \Z(m^+ \wedge n^+) \otimes_{\Z} R(s) \otimes_{\Z} R(t)
  \\
  &\rightarrow \widetilde \Z(m^+ \wedge n^+) \otimes_{\Z} R(s+t) 
  \\ &=
  \widetilde R(m^+ \wedge n^+, s+t), 
  \end{align*}
  induced by the multiplication in $R$ and the unit is given by the
  composition 
  \begin{displaymath}
    \eta : n^+ \rightarrow \widetilde \Z(n^+) \rightarrow
    \widetilde\Z(n^+) \otimes_{\Z} R(0) = \widetilde R(n^+,0)
  \end{displaymath}
  where the last map is induced from the unit of $R$.
\end{example}
The topological Hochschild homology of a filtered FSP $L$ is the filtered
pointed simplicial set $\THH(L)$ with $k$-simplices of $\THH(L)(s)$ given by
the 
homotopy colimit 
$$\hocolim{(i_0,\dots,i_k) \in I^{k+1}} F (S^{i_0} \wedge \dots \wedge S^{i_k},  (L(S^{i_0})
\wedge \dots \wedge L(S^{i_k}))(s)),$$
where the smash product of the $L(S^{i_\alpha})$'s is a smash product of
filtered pointed simplicial sets, and with cyclic structure of
Hochschild type.
We define cyclic spaces $\overline{\THH}(L,s)$
for $s \in \Z$ with $k$-simplices given by the
homotopy colimit
$$\hocolim{(i_0,\dots,i_k) \in I^{k+1}} F (S^{i_0} \wedge \dots \wedge
S^{i_k},  \frac{(L(S^{i_0}) 
\wedge \dots \wedge L(S^{i_k}))(s)}{(L(S^{i_0})
\wedge \dots \wedge L(S^{i_k}))(s-1)}),$$
and with cyclic structure as for the Hochschild construction.

Of course there is also a filtered version of the Gamma space $n^+
\mapsto \THH(L;n^+)$ with $k$-simplices of $\THH(L;n^+)(s)$ given by the
homotopy colimit: 
$$\hocolim{(i_0,\dots,i_k) \in I^{k+1}} F (S^{i_0} \wedge \dots \wedge S^{i_k},  (L(S^{i_0})
\wedge \dots \wedge L(S^{i_k}))(s)\wedge n^+),$$
and there is a Gamma space $n^+ \mapsto \overline \THH(L,s;n^+)$ where
the 
$k$-simplices of $\overline \THH(L,s;n^+)$ are given by the homotopy colimit
$$\hocolim{(i_0,\dots,i_k) \in I^{k+1}} F (S^{i_0} \wedge \dots \wedge S^{i_k},  \frac{(L(S^{i_0})
\wedge \dots \wedge L(S^{i_k}))(s)}{(L(S^{i_0})
\wedge \dots \wedge L(S^{i_k}))(s-1)}\wedge n^+).$$
\begin{lem}
\label{conGamma}
  Let $X$ be a Gamma space satisfying that $X(n^+)$ is $l$-connected
  for every $n\ge 0$. If $U$ is $m$-connected and $V$ is
  $n$-connected, then the map $X(U) \wedge V \rightarrow X(U \wedge
  V)$ is $2m+n+l+3$-connected provided that $l,m$ and $n$ are $> 1$.
\end{lem}
\begin{proof}
  We consider the cofibre $X(U \wedge V) / X(U) \wedge V$ as a
  bisimplicial set $i,j \mapsto Z_{ij} = X(U_i \wedge V_i)_j /
  X(U_i)_j \wedge V_i$. Since the cofibre of a cofibration of
  $l$-connected spaces is $l$-connected $Z_{i\cdot}$ is $l$-connected
  for every $i$, and by lemma \ref{approxlemma} $Z_{\cdot j}$ is
  $2m+n+3$-connected for every $j$. Using the spectral sequence of
  Bousfield-Friedlander \cite[thm. B.5]{Bousfield-Friedlander} we
  obtain the assertion of the lemma. 
\end{proof}
Together with the approximation lemma of B\"okstedt (\cite{Bokstedt1} or
\cite[lemma 2.5.1]{Brun}) the above lemma proves the following.
\begin{lem}
\label{kapprox}
  Given a pointed functor $X : \Gamma^{k+1} \rightarrow \sset$ and
  $(j_0,\dots , j_k) \in I^{k+1}$ the map
  \begin{multline*}
    F(S^{j_0} \wedge \dots \wedge S^{j_k}, X(S^{j_0}, \dots , S^{j_k})) 
   \rightarrow \\
    \hocolim{(i_0,\dots,i_k) \in I^{k+1}} 
    F(S^{i_0} \wedge \dots \wedge S^{i_k}, X(S^{i_0}, \dots , S^{i_k})) 
  \end{multline*}
  is $j-1$-connected, where $j$ denotes the minimum of $j_0,\dots,j_k$.
\end{lem}
\begin{lem}
  The spectra $n \mapsto \THH(L;S^n)(s)$ and $n \mapsto \overline
\THH(L,s;S^n)$ are $\Omega$-spectra.
\end{lem}
\begin{proof}
  Let us note that lemma \ref{conGamma} gives that ${(L(S^{i_0}) 
  \wedge \dots \wedge L(S^{i_k}))(s)}$ is $i_0 + \dots + i_k
  -1$-connected. 
  Replacing lemma \ref{approxlemma} by lemma \ref{kapprox} the proof of
  lemma \ref{thhisomega} also proves this lemma.
\end{proof}
We shall say that a filtered FSP $L$ is filtered by cofibrations if
for every $X$
and $s$
the map $L(X)(s-1) \rightarrow L(X)(s)$ is a cofibration. 
\begin{lem}\label{thhfibreseq}
  Let $L$ be a filtered FSP, filtered by cofibrations. Then the map
  from $\THH(L)(s-1)$
  to the homotopy fibre of the map $q: \THH(L)(s) \rightarrow
  \overline{\THH} (L,s)$ is a weak equivalence. 
\end{lem}
\begin{proof}
  Let us start by showing that 
  the map from the mapping cone of
  the map $\THH_k(L;S^n)(s-1)
  \rightarrow \THH_k(L;S^n)(s)$ to $\overline
  \THH_k(L;S^n,s)$ is $2(n-1)$-connected.
  If $X \rightarrow Y$ is a cofibration of
  $i+n$-connected pointed simplicial sets then by applying the 
  Blakers--Massey  theorem \cite[p. 366]{Whitehead} several times we
  see that the 
  map from the mapping 
  cone of the map $F(S^i,X) \rightarrow F(S^i,Y)$ to $F(S^i,Y/X)$ is
  $2(n-1)$-connected.
  Since the mapping cone
  construction commutes with geometric realization it follows that
  the map from the mapping cone of
  the map $\THH_k(L;S^n)(s-1)
  \rightarrow \THH_k(L;S^n)(s)$ to $\overline
  \THH_k(L;S^n,s)$ is $2(n-1)$-connected.

  Now let $q(S^n)$ denote
  the map $\THH(L;S^n)(s) \rightarrow \overline
  \THH(L,s;S^n)$, and let $hFq(S^n)$ denote its homotopy fibre. It then
  follows from the Blakers--Massey theorem that the map
  $\THH(L;S^n)(s-1) \rightarrow hFq(S^n)$ is $2n-4$-connected. From
  the weak equivalence $F(S^n,hFq(S^n)) \simeq hFq$ it
  follows that the map $\THH(L)(s-1) \rightarrow hFq$ is
  $n-4$-connected. Since $n$ is arbitrary, it follows that this map is
  a weak equivalence.
\end{proof}
\begin{remark}\label{th0/th-1}
Given a filtered FSP $L$ there is an FSP $L(0)/L(-1)$ taking $X$ to
$L(X)(0)/L(X)(-1)$. If $L(s) = L(0)$ when $s \ge 0$ then by
proposition \ref{idquotients} $\overline
\THH(L,0) \cong \THH(L(0)/L(-1))$ and $\THH(L)(0) = \THH(L(0))$. 
In this case the above lemma says that $\THH(L)(-1)$ is weakly
equivalent to the 
homotopy fibre of the map $\THH(L(0)) \rightarrow \THH(L(0)/L(-1))$.  
\end{remark}

\section{Cyclotomic structure}
\label{cyclotomic}
In this section we shall describe how the filtration on topological Hochschild
homology of an FSP filtered by cofibrations is compatible with
topological cyclic 
homology.
We have based our presentation on the elementary version of
topological cyclic homology given in \cite{BHM}.
Alternatively we could use the cyclotomic
spectra in the sense of Madsen \cite{Bokstedt-Madsen}. Since we do not
need them for the main 
result of this paper we 
have chosen the technically simpler version of $TC$.

\subsection{Gamma epicyclic spaces}
\label{gammaepicyclic}
Let us recall Goodwillie's notion of an epicyclic space from \cite{MSRI}.
\begin{defn}
  An epicyclic space is a cyclic space $Y$ equipped with maps $r_q :
  Y_{qj-1}^{C_q} \rightarrow Y_{j-1}$ for all $q \ge 1$ and $j \ge 1$,
  satisfying:
  \begin{enumerate}
  \item $r_q : (\sd_q Y)^{C_q} \rightarrow Y$ is cyclic.
  \item $r_a \circ r_q = r_{aq}: (\sd_{aq} Y)^{C_{aq}} \rightarrow Y$.
  \item $r_1$ is the identity.
  \end{enumerate}
\end{defn}
Here $\sd_q Y$ denotes the $q$-fold edgewise subdivision of $Y$ with
$j$-simplices $(\sd_q Y)_j = Y_{qj-1}$. For a
treatment of edgewise subdivision we refer to \cite{BHM}. The most important
properties of edgewise subdivision are that there is a simplicial
action of $C_q$ on $\sd_q Y$, that there is an action of $S^1$ on
$|\sd_q Y|$ extending the simplicial action of $C_q$, that there is an
$S^1$-isomorphism 
$|\sd_q Y| \cong |Y|$, and that $\sd_{aq} = \sd_a \sd_q$.
Note that $r_q$ induces a $C_a$-equivariant map $(\sd_{aq}Y)^{C_{q}} =
\sd_a(\sd_q Y)^{C_q} \rightarrow \sd_a Y$ for any $a$. 

Write $Y^{C_a}$ for the topological space $|(\sd_a Y)^{C_a}| \cong
|Y|^{C_a}$. There is a map $f_q :Y^{C_{aq}} \cong |Y|^{C_{aq}}
\rightarrow |Y|^{C_a} \cong Y^{C_a}$ induced by inclusion of fixed
points. We shall call this map the {\em Frobenius map}.
The map $(\sd_{aq}Y)^{C_{q}} = \sd_a(\sd_q
Y)^{C_q} \rightarrow \sd_a Y$ induces a map
$r_q : Y^{C_{aq}} = |({(\sd_{aq}Y)^{C_{q}}})^{C_a}|  \rightarrow
|(\sd_a Y)^{C_a}| 
= Y^{C_a}$, and we will call $r_q$ the $q$'th {\em restriction
map}. (Following Hesselholt and 
Madsen \cite{Hesselholt-Madsen}, in conflict with Goodwillie's
terminology) 
The maps $r_{q'}$ and $f_q$ commute, that is $f_q r_{q'} = r_{q'} f_q$.

Let us fix a prime $p$.
The restriction and Frobenius maps induce maps $r,f : \prod_{n\ge 0}
Y^{C_{p^n}} \rightarrow \prod_{n\ge 0} Y^{C_{p^n}}$. We let $tr(Y,p)$
denote the homotopy equalizer of $r$ and the identity. Since $r f = f
r$, the map $f$ induces an endomorphism on $tr(Y,p)$. We define
$tc(Y,p)$ to be the homotopy equalizer of $f$ and the identity on
$tr(Y,p)$. Note that since homotopy limits commute we could equally well have
interchanged the roles of $r$ and $f$ in the definition of $tc(Y,p)$.

\begin{defn}
A Gamma epicyclic space is a Gamma object in the category of
epicyclic spaces. 
\end{defn}

The main example of a Gamma epicyclic space is topological
Hochschild homology. The restriction map $r_q : \sd_q
\THH(L,n^+)^{C_q} \rightarrow \THH(L,n^+)$ is defined degreewise by
the following chain of maps:
\begin{align*}
& (\sd_q \THH(L,n^+))^{C_q}_k \cong \\
&  \hocolim{(n_0,\dots,n_k) \in I^{k+1}} F((S^{n_0} \wedge \dots \wedge S^{n_k})^{\wedge q},
    (L(S^{n_0}) \wedge \dots \wedge L(S^{n_k}))^{\wedge q} \wedge
    n^+)^{C_q} \rightarrow \\
&  \hocolim{(n_0,\dots,n_k) \in I^{k+1}} F(((S^{n_0} \wedge \dots
\wedge S^{n_k})^{\wedge 
    q})^{C_q} ,
    (L(S^{n_0}) \wedge \dots \wedge L(S^{n_k}))^{\wedge q} \wedge
    n^+)^{C_q}) \cong \\
&  \hocolim{(n_0,\dots,n_k) \in I^{k+1}} F(S^{n_0} \wedge \dots \wedge S^{n_k},
    L(S^{n_0}) \wedge \dots \wedge L(S^{n_k}) \wedge
    n^+) = \\
&    \THH(L,n^+)_k.  
\end{align*}
The first isomorphism is due to the isomorphism
$(\hocolim{I^{q(k+1)}} Z)^{C_q} \cong \hocolim{I^{k+1}} Z^{C_q}$. The
second map is given by restriction to fixed points and the last
isomorphism is induced by the point set isomorphism $(X^{\wedge
q})^{C_q} \cong X$.

Given a Gamma epicyclic space $X$, 
we obtain simplicial epicyclic spaces $X(S^n)$. We can view these as
epicyclic spaces and consider the spaces $tc(X(S^n),p)$. In order to
see that these spaces assemble to a spectrum, let us first note that
$n \mapsto (X(S^n))^{C_a} \cong |X(S^n)|^{C_a}$ is a spectrum because
the category of spectra is closed under limits, and limits are
constructed degreewise. Since the same remark applies to homotopy
limits we have a spectrum $n \mapsto tc(X(S^n),p)$. 
\begin{defn}
  Let $X$ be a Gamma epicyclic space. {\emph Topological cyclic
  homology at the prime $p$} of $X$ is the spectrum $TC(X,p)$ with
  $n$'th space 
  $TC(X,p)_n = tc(X(S^n),p)$. 
\end{defn}

We shall write $TC(L,p)$ instead of $TC(\THH(L),p)$.
Our definition of topological cyclic homology at the prime $p$ agrees
with the definition of B\"okstedt, Hsiang and Madsen
\cite[def. 5.12.]{BHM} At this point it is clear that our version
of $TC$ has the same underlying space as the one in \cite{BHM}. To see
that the deloopings agree we first note that our spectrum $TC(L,p)$
is stably equivalent to 
the ones in Goodwillies note \cite{Goodwillie} and in \cite[definiton
4.1]{Hesselholt-Madsen}. Next we appeal to 
\cite[prop. 2.6.2.]{Hesselholt-Madsen}. 

By a Gamma cyclic space we shall mean a Gamma object in the category
of cyclic pointed spaces. Given a Gamma cyclic space $X$ and a closed subgroup
$H$ of $S^1$ we shall let $X^H$ denote the spectrum $n \mapsto
|X(S^n)|^H$. If $X$ is a Gamma epicyclic space the restriction and
Frobenius maps $r_r,f_q : X(S^n)^{C_{aq}} \rightarrow X(S^n)^{C_a}$
induce maps $R_q,F_q : X^{C_{aq}} \rightarrow X^{C_a}$ of
spectra. Given a cyclic pointed space $Z$ we shall define a spectrum
$Z \wedgel X$ by the formula
\begin{displaymath}
  (Z \wedgel X)_n = \colim{W \subset \U} \pMap(|S^{W-\R^n}|, |Z| \wedge |X(S^W)|).
\end{displaymath}
Here $\U$ denotes a complete $S^1$-universe (e.g. $\U = \bigoplus_{n
\in \Z, \alpha \in \N} {\mathbb C}(n)_{\alpha}$) and the colimit runs over
finite dimensional sub inner spaces $W$ of $\U$ containing $\R^n$. The
symbol $W - \R^n$ denotes the orthogonal complement of $\R^n$ in $W$
and $S^W$ denotes the singular complex of the one point
compactification of $W$. There are several possible
actions of $S^1$ on $|X(S^W)|$. Using the functoriality of $X$ the
action of $S^1$ on $S^W$ induces an action of $S^1$ on $X(S^W)$. On
the other hand forgetting the action of $S^1$ on $S^W$ we still have a
cyclic structure on $X(S^W)$ given rise to an action of $S^1$ on
$|X(S^W)|$. The two actions just described commute and therefore we
end up with an action of $S^1 \times S^1$ on $|X(S^W)|$. We shall
always let $S^1$ act on $|X(S^W)|$ by pulling back the action of $S^1
\times S^1$ along the diagonal $S^1 \rightarrow S^1 \times S^1$.
Letting $S^1$ act on the pointed mapping space $\pMap(|S^{W - \R^n}, |Z| \wedge
|X(S^W)|)$ by conjugation we obtain an action of $S^1$ on $(Z \wedgel
X)_n$ and we obtain a spectrum $(Z \wedgel X)^H$ with $(Z \wedgel
X)^H_n = ((Z \wedgel X)_n)^H$ for every closed subgroup $H$ of $S^1$.
Since the map $U \wedge X(V) \rightarrow X(U \wedge V)$ is a
cofibration for every $U$ and $V$ the above construction is
homotopically meaningful. (See the discussion in \cite[Appendix
A]{Hesselholt-Madsen}.) 

In particular we can consider the spectrum $(S^0 \wedgel X)^H$. There
is a map $X^H \rightarrow (S^0 \wedgel X)^H$ induced by the map
$|X(S^n)| \rightarrow \pMap(|S^W|,|X(S^n \wedge S^W)|)$. 
According to \cite[prop. 2.4]{Hesselholt-Madsen} this map is an
equivalence when $X = \THH(L)$ and $H$ is finite.

Using the standard cyclic model of $S^1$ we can consider $ES^1$ as a
cyclic space, and we can consider the spectrum $(ES^1_+ \wedgel
X)^C$. According to \cite[thm. 7.1. p. 97]{LMS} it represents the
$C$-homotopy orbit spectrum $X_{hC}$ of $X$ in the homotopy category
when $C$ is finite and $(ES^1_+ \wedgel X)^{S^1}$ represents the
suspension $S^1 \wedge X_{hS^1}$ of the $S^1$-homotopy orbits of
$X$. From now on we shall always use these representatives for
homotopy orbits. The inclusions of fixed points $(ES^1_+ \wedgel
X)^{C_{aq}} \rightarrow (ES^1_+ \wedgel X)^{C_a}$ and $(ES^1_+ \wedgel
X)^{S^1} \rightarrow (ES^1_+ \wedgel X)^{C_a}$ represent the transfer
maps $\trf_q : X_{hC_{qa}} \rightarrow X_{hC_a}$ and $\trf_\infty :
S^1 \wedge X_{hS^1} \rightarrow X_{hC_a}$ respectively. 
We shall always use these representatives for the the transfer maps.

\begin{defn}
\label{defpcyc}
  A $p$-cyclotomic Gamma space is a Gamma epicyclic space $X$ satisfying
  the following two conditions.
  \begin{enumerate}
    \item  The map $X^C \rightarrow (S^0 \wedgel X)^C$ is an
    equivalence for every finite $p$-subgroup $C$ of $S^1$.
    \item  The norm map $N : X_{hC_{p^n}} \rightarrow X^{C_{p^n}}$
  defined as the 
  composite 
  \begin{displaymath}
    X_{hC_{p^n}} = (ES^1_+ \wedgel X)^{C_{p^n}} \rightarrow (S^0
  \wedgel X)^{C_{p^n}}  \simeq X^{C_{p^n}}
  \end{displaymath}
  fits into a cofibration sequence
$    X_{hC_{p^n}} \overset{N} \rightarrow X^{C_{p^{n}}} \overset{R_p}
    \rightarrow   X^{C_{p^{n-1}}} $
  for every $n \ge 1$.
  \end{enumerate}
\end{defn}

Note that the norm map is only defined in the homotopy category and
that the diagram
\begin{displaymath}
  \begin{CD}
    X_{hC_{p^n}} @>N>> X^{C_{p^n}} \\
    @V{\trf_p}VV @VV{F_p}V \\
    X_{hC_{p^{n-1}}} @>N>> X^{C_{p^{n-1}}} 
  \end{CD}
\end{displaymath}
commutes. It is proven in \cite[lemma 2.5]{Hesselholt-Madsen} and
\cite[prop. 2.4]{Hesselholt-Madsen} that $\THH(L)$
satisfies $(1)$ and $(2)$ above. In conclusion $\THH(L)$ is a $p$-cyclotomic
Gamma space.

Below we shall use the following
lemma due to Goodwillie. It can be found in \cite{Madsentrc} as lemma
4.4.9. 
\begin{lem}
\label{S1transfer}
  For any epicyclic Gamma space $X$ the $S^1$-transfer induces a map
  \begin{displaymath}
    S^1 \wedge X_{hS^1} \rightarrow \holim{\trf_p} X_{hC_{p^n}}.
  \end{displaymath}
  This map becomes an equivalence after $p$-completion.
\end{lem}
Let us sketch an alternative proof of this lemma. Since $(ES^1_+)_k =
(S^1_+)^{\wedge k+1}$ it suffices to show that the map 
$  ((S^1_+)^{\wedge k+1} \wedgel X)^{S^1} \rightarrow \holim{F_p}
  ((S^1_+)^{\wedge k+1} \wedgel X)^{C_{p^{n}}}$
is an equivalence for every $k \ge 0$. There is an isomorphism
$(S^1_+)^{\wedge k+1} \wedgel X \cong S^1_+ \wedgel ((S^1_+)^{\wedge
k} \wedge X)$, where $(S^1_+)^{\wedge
k} \wedge X$ denotes the Gamma cyclic space $n^+ \mapsto (S^1_+)^{\wedge
k} \wedge X(n^+)$. Therefore the proof of lemma \ref{S1transfer}
reduces to
showing that the map $(S^1_+ \wedgel X)^{S^1} \rightarrow \holim{F_p}
  (S^1_+ \wedgel X)^{C_{p^n}}$ becomes an equivalence after completion
at $p$. This is the statement of  
\cite[lemma 8.2]{Hesselholt-Madsen}.

\subsection{Cyclotomically filtered Gamma spaces}
In this section we shall present a filtered version of $p$-cyclotomic
Gamma spaces. Let us begin with a filtered version
of the notion of an epicyclic space.

\begin{defn}
  An epicyclic filtered space is a filtered cyclic space $Y$ equipped
  with maps $r_q : 
  Y_{qj-1}(s)^{C_q} \rightarrow Y_{j-1}([s/q])$ for all $q,j \ge 1$ and
  $s \in \Z$, 
  satisfying:
  \begin{enumerate}
  \item $r_q : (\sd_q Y)(s)^{C_q} \rightarrow Y([s/q])$ is cyclic.
  \item $r_a \circ r_q = r_{aq}: (\sd_{aq} Y)(s)^{C_{aq}} \rightarrow
    Y([s/(aq)])$. 
  \item $r_1$ is the identity.
  \end{enumerate}
\end{defn}

Here $[s/q]$ denotes the greatest integer $\le s/q$.  Write
$Y^{C_a}(s)$ for the topological space $|(\sd_a Y(s))^{C_a}| \cong
|Y(s)|^{C_a}$. There is a Frobenius map $f_q :Y^{C_{aq}}(s) \cong
|Y(s)|^{C_{aq}} \rightarrow |Y(s)|^{C_a} \cong Y^{C_a}(s)$ induced by
inclusion of fixed points.

A Gamma epicyclic filtered space is a Gamma object in the category of
epicyclic filtered spaces.

Topological
Hochschild homology of an FSP filtered by cofibrations is the main
example of a Gamma  epicyclic filtered  space. The restriction 
map 
$$r_q : \sd_q 
\THH(L,n^+)(s)^{C_q} \rightarrow \THH(L,n^+)([s/q])$$ 
is defined degreewise by
the following chain of maps:
\begin{align*}
&  (\sd_q \THH(L,n^+)(s))^{C_q}_k \cong \\
&  \hocolim{(n_0,\dots,n_k) \in I^{k+1}} F((S^{n_0} \wedge \dots \wedge S^{n_k})^{\wedge q},
    (L(S^{n_0}) \wedge \dots \wedge L(S^{n_k}))^{\wedge q}(s) \wedge
    n^+)^{C_q} \rightarrow \\
&  \hocolim{(n_0,\dots,n_k) \in I^{k+1}} F(((S^{n_0} \wedge \dots \wedge S^{n_k})^{\wedge
    q})^{C_q} ,
    (L(S^{n_0}) \wedge \dots \wedge L(S^{n_k}))^{\wedge q}(s) \wedge
    n^+)^{C_q}) \\
&  \cong \hocolim{(n_0,\dots,n_k) \in I^{k+1}} F(S^{n_0} \wedge \dots \wedge S^{n_k},
    (L(S^{n_0}) \wedge \dots \wedge L(S^{n_k}))([s/q]) \wedge
    n^+) = \\
    & \THH(L,n^+)([s/q])_k.  
\end{align*}
The first isomorphism is due to the isomorphism
$(\hocolim{I^{q(k+1)}} Z)^{C_q} \cong \hocolim{I^{k+1}} Z^{C_q}.$ 
The
second map is given by restriction to fixed points and the last
isomorphism is induced by the point set isomorphism $(X^{\wedge
q})(s)^{C_q} \cong X([s/q])$ for $X$ a space filtered by cofibrations.
This last isomorphism is not obvious though, so we state it
as a lemma.
\begin{lem}
\label{filtfix}
  Let $Y$ be a filtered space, filtered by cofibrations. There is an
  isomorphism $\left( Y^{\wedge 
  q}(s) \right)^{C_q} \cong Y([s/q])$
\end{lem}
\begin{proof}
  The diagonal induces a map
\begin{displaymath}
  Y([s/q]) \overset{\cong}\rightarrow \left( (Y^{\wedge q})(s) \right)^{C_q}.
\end{displaymath}
We will show that this is an isomorphism of simplicial sets. We may
  assume that $Y$ is a discrete set filtered by injections.
We note that by the pushout diagram in the proof
of lemma \ref{idquotients} the map
$(Y^{\wedge q})(i) \rightarrow (Y^{\wedge q})(i+1)$ is an injection
for all $i \in \Z$, and therefore we have an injection
\begin{displaymath}
  Y^{\wedge q}(s) \hookrightarrow Y^{\wedge q}(\infty) \cong
  (Y(\infty))^{\wedge q},
\end{displaymath}
with the convention that $Y^{\wedge q}(\infty) = \colim{i} Y^{\wedge
q}(i)$ and $Y(\infty) = \colim{i} Y(i)$. There is a commutative
diagram
  $$\begin{matrix} 
    Y([s/q]) 
  &\to& 
    (Y^{\wedge q}(s))^{C_q} 
  \\
  \downarrow && \downarrow \\
    Y(\infty) 
  & \to & 
    (Y(\infty)^{\wedge q})^{C_q}
  \end{matrix}$$
where the vertical arrows are injections. It follows from the diagram
that the map $Y([s/q]) \rightarrow (Y^{\wedge q}(s))^{C_q}$ is
injective. To see that it is onto, let us pick a representative
$((a_1,\dots,a_q),(y_1,\dots,y_q))$ for a point $y$ in 
\begin{displaymath}
  (Y^{\wedge q}(s)) = \colim{a_1+\dots+a_q \le s}Y(a_1) \wedge \dots
  \wedge Y(a_q),
\end{displaymath}
fixed under the
$C_q$-action. From the condition $a_1 + \dots + a_q \le s$, it follows
that there exists an $i$ such that $a_i \le s/q$. Since the image of
$y$ in $Y(\infty)^{\wedge q}$ is a fixed point, we must have that
$(a_1,y_1), \dots , (a_q,y_q)$ represent the same element in
$Y(\infty)$. Since the map $Y^{\wedge q}(s) \rightarrow Y^{\wedge
  q}(\infty)$ is injective, it follows that
$((a_i,\dots,a_i),(y_i,\dots,y_i))$ represents $y$, and we can conclude
that the map $Y([s/q]) \rightarrow (Y^{\wedge q}(s))^{C_q}$ is onto.
\end{proof}
\begin{defn}
  A $p$-cyclotomic filtered Gamma space is a Gamma epicyclic filtered
  space $X$ satisfying 
  the following two conditions.
  \begin{enumerate}
    \item  The map $X^C(s) \rightarrow (S^0 \wedgel X(s))^C$ is an
    equivalence for every finite $p$-subgroup $C$ of $S^1$ and 
  $s \in \Z$. 
    \item  The norm map $N : X(s)_{hC_{p^n}} \rightarrow X^{C_{p^n}}(s)$
  defined as the composite 
  \begin{displaymath}
    X(s)_{hC_{p^n}} = (ES^1_+ \wedgel X(s))^{C_{p^n}} \rightarrow (S^0
  \wedgel X(s))^{C_{p^n}}  \simeq X^{C_{p^n}}(s)
  \end{displaymath}
  fits into a cofibration sequence
$    X(s)_{hC_{p^n}} \overset{N} \rightarrow X^{C_{p^{n}}}(s) \overset{R_p}
    \rightarrow   X^{C_{p^{n-1}}}([s/p]) $
  for every $n \ge 1$.
  \end{enumerate}
\end{defn}
The proof of \cite[prop 2.4]{Hesselholt-Madsen} shows that for any filtered FSP
$L$ the Gamma cyclic space  $\THH(L)(s)$ satisfies $(1)$ above. If $L$
is filtered by cofibrations we 
have by lemma \ref{filtfix} that 
$$((L(S^{i_0}) \wedge \dots \wedge L(S^{i_k}))^{\wedge r}(s))^{C_q}
\cong  (L(S^{i_0}) \wedge \dots \wedge L(S^{i_k}))^{\wedge r/q}([s/q]).$$
It follows from lemma \ref{conGamma} that its connectivity is at least
$(i_0+\dots +i_k)r/q-1$. The proof of \cite[prop
2.5]{Hesselholt-Madsen} together 
with the above observation shows that $\THH(L)$ satisfies $(2)$ in the
above definition. Hence $\THH(L)$ is a $p$-cyclotomic filtered Gamma
space.

\section{Filtered topological cyclic homology}

Given a Gamma epicyclic filtered space $X$, 
the Gamma spaces $X(-1)$, $X(0)$ and $X(\infty) = \colim{s} X(s)$ come
equipped with an epicyclic structure. In this section we shall define
a filtered Gamma space $TC = TC(X,p)$, the topological cyclic homology of
$X$. This will be a generalization of non-filtered topological cyclic
homology in the sense that $TC(s)$ is
the (non-filtered) topological cyclic homology $TC(X(s),p)$ of the
Gamma epicyclic space $X(s)$ for $s = 
-1, 0, \infty$. Here we use the notation $TC(\infty) = \colim{s}
TC(s)$. The sign of $s$ plays an important role in the definition of
$TC(s)$. In fact $TC(s)$ is defined using only the $X(r)$ where $r$
has the same sign as $s$. 

We start by defining $tc(Y)(s)$ for a filtered epicyclic space $Y$.
Let us start by defining $tc(Y)(s)$ when $s \le 0$. In this case $s \le
[s/p]$. (Recall that $[s/p]$ is the greatest integer less that or
equal to $s/p$.) 
Therefore there is an inclusion $Y(s) \overset i
\rightarrow Y([s/p])$. Using this inclusion and the restriction map
$Y(s)^{C_{p^n}} \overset{r_p} \rightarrow Y([s/p])^{C_{p^{n-1}}}$ we
obtain maps
\begin{displaymath}
  \prod_{n \ge 0} Y(s)^{C_{p^n}} \overset{\textstyle \overset{r(s)}
  \rightarrow} 
  {\underset{i(s)} \rightarrow} \prod_{n \ge 0} Y([s/p])^{C_{p^n}}.
\end{displaymath}
Let us define $tr(Y)(s)$ to be the homotopy equalizer of $r(s)$ and
$i(s)$. The Frobenius maps 
$Y(s)^{C_{p^n}} \overset {f(s)} \rightarrow Y(s)^{C_{p^{n-1}}}$ commute
with $r(s)$ and $i(s)$, and therefore they induce an endomorphism $f(s)$ of
$tr(Y)(s)$. We define $tc(Y)(s)$ to be the homotopy equalizer of $f(s)$ and
the identity. When $s = -1,0$ this definition of $tc(Y)(s)$ agrees with
the definition of $tc(Y(s))$ given in section \ref{gammaepicyclic}.
There is an alternative definition of $tc(Y)(s)$ where we interchange
the roles of $r$ and $f$ going as follows: We let
$tf(Y)(s)$ denote the homotopy equalizer of the maps
\begin{displaymath}
  \prod_{n \ge 0} Y(s)^{C_{p^n}} \overset{\textstyle \overset{f(s)}
  \rightarrow} 
  {\underset{\id} \rightarrow} \prod_{n\ge 0} Y(s)^{C_{p^n}}.
\end{displaymath}
$tc(Y)(s)$ is equivalent to the homotopy equalizer of the
maps $r(s)$ and $i(s)$ from $tf(Y)(s)$ to $tf(Y)([s/p])$. 

Now let us define $tc(Y)(s)$ when $s \ge 0$. The restriction map induces
an endomorphism $r(s)$ on the product $\prod_{n\ge 0}
Y(sp^n)^{C_{p^n}}$. We define 
$tr(Y)(s)$ to be the homotopy equalizer of $r(s)$ and the identity. The
Frobenius map induces a map $f(s) : tr(Y)(s) \rightarrow tr(Y)(sp)$. Since
$s \le sp$, there is an inclusion $Y(s) \rightarrow Y(sp)$ inducing a
map $i(s) : tr(Y)(s) \rightarrow tr(Y)(sp)$. We define $tc(Y)(s)$ to be the
homotopy equalizer of $f(s)$ and $i(s)$. When $s = 0$ this definition
of $tc(Y)(s)$ agrees with the definition given above.

Given a Gamma epicyclic filtered space $X$
we have spaces 
$tc(X(S^n))(s)$.
Filtered topological cyclic homology at the prime $p$ of $X$ is the spectrum
$TC = TC(X,p)$ with $TC(X,p)(s)_n = tc(X(S^n),p)(s)$.
Similarly let $TR = TR(X,p)$ be the filtered spectrum with 
$TR(X,p)(s)_n
= tr(X(S^n),p)(s).$
Let us note that $TC(X(s),p) \cong TC(X,p)(s)$ for $s =
-1,0,\infty$. This fact together with the two following lemmas is our
justification for the definition of $TC(X,p)$. For $s = -1,0,\infty$
the spectrum $TR(X(s),p)$ can be 
rewritten as the sequential homotopy 
limit of $X(s)^{C_{p}}$  with respect to the restriction maps.
There is no such rewriting possible for $s \ne -1,0,\infty$.

\begin{lem}
  Let $X$ be a $p$-cyclotomic filtered Gamma space filtered by cofibrations. 
  Suppose that the connectivity of the map $X(s) \rightarrow
  X(\infty)$ tends to infinity as $s$ grows. Then $X(\infty)$ is a
  $p$-cyclotomic Gamma space and  $TC(X,p)(\infty)$ is
  stably equivalent to $TC(X(\infty),p)$.
\end{lem}
Recall that $X(\infty)$ is the Gamma epicyclic space space with underlying
Gamma space $\colim{s} X(s)$, and that $TC(\infty) = \colim{s}
TC(s)$.
\begin{proof}
  Suppose that the map $X(s) \rightarrow X(\infty)$ is $k$-connected
  for $s \ge N \ge 0$. Using the cofibration sequence
  \begin{displaymath}
    X(sp^n)_{hC_{p^n}} \rightarrow X(sp^n)^{C_{p^n}} \rightarrow
  X(sp^{n-1})^{C_{p^{n-1}}} 
  \end{displaymath}
  we can by induction show that the map $X(sp^n)^{C_{p^{n}}}
  \rightarrow X(\infty)^{C_{p^n}}$ is $k$-connected when $s \ge N$,
  for all $n$. It follows that the map
  \begin{displaymath}
    \prod_{n \ge 0} X(sp^n)^{C_{p^n}} \rightarrow \prod_{n \ge 0}
  X(\infty)^{C_{p^n}} 
  \end{displaymath}
  is $k$-connected when $s \ge N$. Therefore the map $TR(s)
  \rightarrow TR(X(\infty),p)$ is at least $k-1$-connected, and the map
  $TC(s) \rightarrow TC(X(\infty),p)$ is at least $k-2$-connected when
  $s \ge N$.
\end{proof}
The next lemma says that if $X$ is a $p$-cyclotomic filtered Gamma
space and the connectivity of $X(s)$ tends to
infinity as $s$ decreases then $\holim{s} TC(s)$ is 
contractible. 
\begin{lem}
  \label{lemtrtc}
  Let $X$ be a $p$-cyclotomic filtered Gamma space. If $s \le 0$ then
  $TR(s)$ and $\Sigma TC(s)$ are at least as highly connected as $X(s)$. 
\end{lem}
\begin{proof}
  Let $TR^m(s)$ denote the homotopy equalizer of the diagram
  \begin{displaymath}
    \prod_{0 \le n \le m} X(s)^{C_{p^n}} 
    \overset{\textstyle \overset{R(s)} \rightarrow} {\underset{I(s)}
      \rightarrow}  
    \prod_{0 \le n \le m-1} X([s/p])^{C_{p^n}},
  \end{displaymath}
  where $I(s)$ forgets the $m$'th coordinate, and otherwise $I(s)$ and
  $R(s)$ are truncations of the maps defining $TR(s)$.
  There is an obvious map $TR^m(s) \rightarrow
  TR^{m-1}(s)$
  induced by
  projection away from the last factors of the products. Using the
  norm cofibration sequence the fibre of
  this map may be identified with $X(s)_{hC_{p^m}}$.
  Since homotopy limits commute we have that $TR(s)$ is the
  homotopy limit of the sequence
  \begin{displaymath}
    \dots \rightarrow TR^m(s) \rightarrow TR^{m-1}(s) \rightarrow \dots
    \rightarrow TR^0(s).
  \end{displaymath}
  Since homotopy orbits preserve
  connectivity and $TR^0(s) = X(s)$, $TR(s)$ is a
  sequential homotopy limit of  
  spaces as least as connected as $X(s)$ and with homotopy fibres as
  least as connected as $X(s)$. It follows that
  $TR(s)$ is at least as highly connected as $X(s)$. 
  Since homotopy equalizers at most lower connectivity by one we
  have that $\Sigma TC(s)$ is at
  least as highly connected as $TR(s)$. 
\end{proof}

Given a map $A \rightarrow B$ of spectra we shall denote the
homotopy cofibre by $B/A$. 
\begin{lem}
\label{smartquotients}
  Let $X$ be a $p$-cyclotomic filtered Gamma space,
  let $s < 0$, and assume that $ps \le t < s$. After $p$-completion 
  $TC(s) / TC(t)$ is equivalent to $S^1 \wedge (X(s)/X(t))_{hS^1}$.
\end{lem}
\begin{proof}
  Since $s \le [t/p]$, the inclusion $X(s) \rightarrow X([s/p])$
  induces the trivial map from $X(s)/X(t)$ to
  $X([s/p])/X([t/p])$. Therefore the homotopy equalizer of the
  maps induced by $I(s)$ and $R(s)$ on the quotients of the products
  in the definition of $TR$ agrees with the homotopy fibre of the map
  induced by $R(s)$. Using the norm cofibration sequence we can
  identify this fibre with $\prod_{n \ge 0} (X(s)/X(t))_{hC_{p^n}}$.
  Since $F_p N = N \trf_p$, where $N :
  X(s)_{hC_{p^n}} \rightarrow X(s)^{C_{p^n}}$ denotes the norm map,
  and where $\trf_p$ denotes the transfer map, we have
  that $TC(s)/TC(t) \simeq \holim{\trf_p} (X(s)/X(t))_{hC_{p^n}}$. The
  lemma now follows from lemma \ref{S1transfer}.
\end{proof}
Note that the above lemma applies to the filtration
quotients $TC(s)/TC(s-1)$.
For our main theorem the case $s = -1$ and $t = -p$ is of particular
interest. In that case we get by lemma \ref{lemtrtc} that if 
$X(-p)$ is $k$-connected, then after $p$-adic
completion the map 
$$TC(-1) \rightarrow \frac{TC(-1)}{TC(-p)} \simeq S^1 \wedge \left(
  \frac{X(-1)}{X(-p)} 
\right)_{hS^1}$$
is $k$-connected. 

It follows from remark \ref{th0/th-1} that for an FSP $L$ filtered by
cofibrations and with $L(s) = L(0)$ for $s \ge 0$ we have that
$TC(0)/TC(-1)$ is isomorphic to $TC(L(0)/L(-1),p)$.

\section{Relative $K$-theory of nilpotent ideals}
\label{goodwillietype}
In this section we shall prove the following theorem relating relative
$K$-theory and relative cyclic Homology. One good reference for cyclic
homology is the book of Loday \cite{Loday}. 
\begin{thm}
\label{mainthm} 
  Let $R$ be a simplicial ring with an ideal $I$ satisfying $I^m =
  0$. Suppose that $R$ and $R/I$ are flat as modules over $\Z$. Then
  there is an isomorphism of homotopy
  groups of $p$-adic completions
  \begin{displaymath}
    \pi_i K(R,I)^{\wedge}_p \cong \pi_{i-1} HC(R,I)^{\wedge}_p
  \end{displaymath}
  when $0 \le i < p/(m-1) -2$ and a surjection
  \begin{displaymath}
    \pi_i K(R,I)^{\wedge}_p \rightarrow \pi_{i-1} HC(R,I)^{\wedge}_p 
  \end{displaymath}
  when $i < p/(m-1) -1$.  
\end{thm}
Recall that $K(R,I)$ is the homotopy fibre of the map $K(R)
\rightarrow K(R/I)$ and that $HC(R,I)$ is the homotopy fibre of the
map $HC(R) \rightarrow HC(R/I)$.

The proof uses the results of the previous section plus a number of
results about $TC$ proven elsewhere. We shall 
collect the statements of these results for the convenience of the
reader. 
The following result is due to McCarthy \cite{McCarthy}.
\begin{thm}
\label{relativetc}
Suppose $f\colon R\to S$ is a
homomorphism of simplicial rings and that $\pi \sb 0(f)$ is surjective
and has nilpotent kernel. Then the 
diagram 
$$\begin{matrix} K(R)&\to& TC(\widetilde R,p)\\\downarrow &&\downarrow
\\K(S)&\to& 
TC(\widetilde S,p)\end{matrix}$$ is 
homotopy Cartesian after $p$-adic completion.  
\end{thm}
Suppose that $f$ in the above theorem is degreewise surjective and let
$I \subseteq R$ denote its kernel. Let $TC(\widetilde
R,\widetilde I,p)$ denote
the homotopy fibre of the map $TC(\widetilde R,p) \rightarrow
TC(\widetilde S,p)$. Then the theorem
says that the map $K(R,I) \rightarrow TC(\widetilde R, \widetilde
I,p)$ is an equivalence 
after $p$-adic completion. The theorem in particular applies in the
situation where $I$ is a nilpotent ideal in $R$. 
\begin{lem}
\label{compareTandHH} 
  Let $R$ be a ring which is flat as a module over
  $\Z$. Then the map 
  $\pi_i \THH(\widetilde R)^{\wedge}_p \rightarrow \pi_i \HH(R)^{\wedge}_p$ is
  an isomorphism when $i \le 2p-2$.
\end{lem}
\begin{proof}
  In \cite[thm 4.1]{Pirashvili-Waldhausen} Pirashvili and Waldhausen 
  have established a spectral sequence with $E^2$-term 
  $E^2_{s,t} = \HH_s(R, \pi_t \THH(\widetilde \Z,\widetilde R))$
  converging towards $\pi_{s+t}
  \THH(\widetilde R)$. The lemma follows from the fact that $\pi_0
  \THH(\widetilde \Z,\widetilde R) = R$ 
  and 
  that $\pi_i \THH(\widetilde \Z,\widetilde R)^{\wedge}_p = 0$ when $1
  \le i \le 2p-2$. 
\end{proof}
The following result is dual to a result of Cohen and Jones 
\cite[lemma 1.3]{Cohen-Jones}. We give an alternative proof inspired by a
more 
elementary proof due to B\"okstedt. 
\begin{prop}
\label{folklore}
  Let $A$ be a a cyclic object in the
  category of abelian groups, and let $\widetilde A$ denote the
  Gamma cyclic space with $\widetilde A(n^+) = A \otimes_{\Z}
  \widetilde \Z(n^+)$. Then there is a natural
  isomorphism $HC_*(A) \cong \pi_*({\widetilde A}_{hS^1})$.
\end{prop}
Recall that $\widetilde A_{hS^1}$ is the $S^1$ homotopy orbit spectrum
associated to the spectrum $n \mapsto |\widetilde A(S^n)|$.
\begin{proof}
  We refer to \cite[section 6.2]{Loday} for the notation used in this
  proof. 
  It is well known (see e.g. \cite[theorem 6.2.8]{Loday}, or use the
  argument below) that there is
  an isomorphism $HC_*(A) \cong \Tor^{\Lambda^{\op}}_*(\Z,A)$. To see
  that $\pi_*({\widetilde A}_{hS^1})$ also is isomorphic to this
  $\Tor$ group, it suffices to check the usual properties determining
  $\Tor$
  groups up to isomorphism (see e.g. \cite[theorem
  V.6.1]{Cartan-Eilenberg}). Firstly we note that there are 
  isomorphisms 
  $\pi_0({\widetilde A}_{hS^1}) \cong \pi_0(|A|) \cong \Z
  \otimes_{\Lambda^{\op}} A$. Secondly the representable functors
  $\Z\Lambda^n$ 
  with $\Z \Lambda^n([m]) = \Z[\Lambda([m],[n])]$ form a set of
  projective generators for the category of cyclic objects in the
  category of abelian
  groups, that is, every projective object in this category is a
  quotient of a sum of 
  objects of the form $\Z \Lambda^n$. Since $|\Z\Lambda^n| \cong
  |\Z[S^1 \times \Delta^n]|$ we 
  can use lemma \ref{gammaconnected} to see that ${\widetilde
  {\Z\Lambda^n}}_{hS^1} 
  \simeq {\widetilde \Z}$, and hence 
  $\pi_i({\widetilde {\Z\Lambda^n}}_{hS^1}) = 0$ for $i>0$. Thirdly,
  given a short exact sequence $A' \rightarrow A \rightarrow A''$ of
  cyclic objects in the category of abelian groups, we obtain a
  cofibration sequence ${\widetilde 
  A'} \rightarrow {\widetilde A} \rightarrow {\widetilde A''}$ of
  spectra with an action of $S^1$. Since
  homotopy orbits take cofibration sequences to cofibration sequences
  we obtain a long exact sequence of homotopy groups of homotopy orbits.   
\end{proof}
The proposition in particular says that there is an
isomorphism $\pi_i {\widetilde{\HH(R)}}_{hS^1} \cong HC_i(R)$.
To prove theorem \ref{mainthm} we also need the following lemma. 
\begin{lem}
  Let $L$ be an FSP filtered by cofibrations. Suppose that $L(s) =
  L(0)$ for $s \ge 0$, and that there
  exists $m \ge 0$ such that $L(-m) = *$. Then $\pi_k
  \THH(L)(s) = 0$ when $k 
  < -s/(m-1) -1$. 
\end{lem}
\begin{proof}
  Recall that $(L(S^{n_0}) \wedge \dots \wedge L(S^{n_k}))(s)$ is the
  colimit running over $i_0 + \dots + i_k \le s$ of 
  \begin{displaymath}
    L(S^{n_0})(i_0) \wedge \dots \wedge L(S^{n_k})(i_k).
  \end{displaymath}
  If $i_0 + \dots + i_k \le s$ then there exists an $\alpha$ such that
  $i_\alpha \le s/(k+1)$. Therefore the smash product is zero if
  $s/(k+1) < -m+1$, or equivalently if $k < -s/(m-1)-1$, and hence
  $\THH(L)(s)_k = 0$ if $k < -s/(m-1)-1$.
\end{proof}
Using the above lemma and lemma \ref{smartquotients}, or rather the remark
after it, we obtain the following.
\begin{prop}
  Let $L$ be an FSP filtered by cofibrations, and suppose that $L(s) =
  L(0)$ for $s \ge 0$, and that there
  exists $m \ge 0$ such that $L(X)(-m) = *$ for all $X$. Let $TC =
  TC(L,p)$, and $\THH = \THH(L)$. After
  completion at $p$ there are maps
  \begin{displaymath}
    TC(-1) \rightarrow \frac{TC(-1)}{TC(-p)} \simeq 
    S^1 \wedge \left(\frac{\THH(-1)}{\THH(-p)}\right)_{hS^1}
    \leftarrow S^1 \wedge \THH(-1)_{hS^1}. 
  \end{displaymath}
  Here $\THH(-p)$ is $p/(m-1)-2$-connected and $TC(-p)$ is
  $p/(m-1)-3$-connected, and therefore the map pointing to the right is
  $p/(m-1)-2$-connected and the map pointing to the left is
  $p/(m-1)-3$-connected.  
\end{prop}
Now let $R$ denote a ring with an ideal $I$ satisfying that $I^m =
0$. Considering the $I$-adic filtration $0 =I^m \subseteq \dots
\subseteq I \subseteq R$ of $R$ we obtain a filtered ring $R$ with
$R(s) = I^{-s}$ for $s<0$ and with $R(s) = R$ for $s \ge 0$. 
By the construction in \ref{filtringfiltfsp} we obtain a filtered FSP
$\widetilde R$ with $\widetilde R(n^+,s) = \widetilde \Z(n^+)
\otimes_{\Z} R(s)$. As remarked after \ref{thhfibreseq}
$\THH(R)(-1)$ is equivalent to the
homotopy fibre of the map $\THH(R) \rightarrow
\THH({R/I})$. Since homotopy limits commute we have that 
\begin{displaymath}
  TC(R,I,p) = TC(R,p)(-1).
\end{displaymath}
Applying McCarthy's theorem \ref{relativetc} and the above proposition
we obtain an isomorphism
\begin{displaymath}
  \pi_i K(R,I)^{\wedge}_p  \cong \pi_{i-1} (\THH(
  R)(-1)_{hS^1})^{\wedge}_p 
\end{displaymath}
when $i < p/(m-1) -2$ and we have a surjection
\begin{displaymath}
  \pi_i K(R,I)^{\wedge}_p  \rightarrow \pi_{i-1} (\THH(
  R)(-1)_{hS^1})^{\wedge}_p 
\end{displaymath}
when $i < p/(m-1) -1$. 
Using lemma \ref{folklore} and lemma \ref{compareTandHH} 
we can complete the proof of theorem
\ref{mainthm}. 
\begin{proof}[Proof of theorem \ref{mainthm}]
  Let us write $\THH(R,I)$ instead of $\THH(R)(-1)$.
  We have an isomorphism $\pi_iK(R,I)^{\wedge}_p \cong
  \pi_{i-1} (\THH(R,I)_{hS^1})^{\wedge}_p$ when $i <
  p/(m-1)-2$
  and a surjection $\pi_iK(R,I)^{\wedge}_p \rightarrow
  \pi_{i-1} (\THH(R,I)_{hS^1})^{\wedge}_p$ when $i <
  p/(m-1)-1$. 
  By lemma 
  \ref{compareTandHH} there is for $i \le 2p-1$ an isomorphism $\pi_{i-1}
  (\THH(R,I)_{hS^1})^{\wedge}_p 
  \cong \pi_{i-1} (\HH(R,I)_{hS^1})^{\wedge}_p$.
  Splicing these maps with the isomorphism of lemma
  \ref{folklore} we obtain the asserted isomorphism
  $\pi_iK(R,I)^{\wedge}_p \cong  
  \pi_{i-1}HC(R,I)^{\wedge}_p$ when $i < p/(m-1)-2$ and the surjection
  $\pi_iK(R,I)^{\wedge}_p \rightarrow  
  \pi_{i-1}HC(R,I)^{\wedge}_p$ when $i < p/(m-1)-2$. 
\end{proof}

\section{Computations in cyclic homology}
\label{derivedcyclic}
In this section we shall compute some derived cyclic homology groups
of the ring $\Z/p^n$. The definition of derived cyclic homology depends
on the following lemma. A proof can for example be found in
\cite{Brun}.
\begin{lem}
  Let $A$ be a simplicial ring. There exists a weak equivalence $R
  \overset{\simeq} \rightarrow A$ of simplicial rings, where $R$ is
  degreewise free as an abelian group. If $R'$ is another ring with
  underlying degreewise free abelian group, and with a weak
  equivalence $R' \overset{\simeq} \rightarrow A$, then there is a
  chain of weak equivalences between $R$ and $R'$ through simplicial
  rings with underlying degreewise free abelian groups.
\end{lem}
Let $A$ be a simplicial ring, and choose a weak equivalence $R
\overset{\simeq} \rightarrow A$ as in the above lemma. That is, with
$R$ degreewise free as an abelian group. By functoriality of the
Hochschild construction there is a map $\HH(R) \rightarrow \HH(A)$. By
definition $\HH(R)$ is the {\em derived Hochschild homology} of $A$.
(Some authors call it Shukla homology.) By
the above lemma it is unique up to
weak equivalence. We shall call 
$HC(R)$ the {\em derived cyclic homology} of $A$. We shall use the
notation $\widetilde{HC} (A)$ for $HC(R)$.

Given a discrete ring $A$, we can consider it as a constant simplicial
ring. This way we obtain derived cyclic homology of discrete rings.
\begin{prop}
\label{derivedcycliccom}
  For $0 \le i < 2p$ the derived cyclic homology of $\Z/p^n$ is
  given as follows:
  \begin{displaymath}
    \widetilde {HC}_i (\Z/p^n) =
      \begin{cases}
        \Z/p^{nj} & \text{if $i = 2(j-1) <2p$} \\
        0 & \text{if $i<2p$ is odd}
      \end{cases}
  \end{displaymath}
  and the relative cyclic homology groups are:
  \begin{displaymath}
    \widetilde {HC}_i (\Z/p^n,p^{n-1} \Z/p^n) =
      \begin{cases}
        \Z/p^{j} & \text{if $i = 2(j-1) <2p$} \\
        0 & \text{if $i<2p$ is odd.}
      \end{cases}
  \end{displaymath}
\end{prop}
\begin{proof}
  Let us consider $\Delta^1 = \Delta(-,[1])$ as a pointed simplicial
  monoid as follows. Given $\alpha, \beta: [k] \rightarrow [1]$, we
  let $(\alpha \cdot \beta) (j) = \alpha(j) \cdot \beta(j)$. The
  constant map with value $0 \in [1]$ is the base point. There is a
  pointed submonoid $S^0$ of $\Delta^1$ consisting of the constant
  maps. We shall let $R$ denote the subring $\widetilde \Z(S^0)\oplus_{p^n
  \widetilde \Z(S^0)} p^n\widetilde \Z(\Delta^1)$ of the pointed monoid ring
  $\widetilde \Z(\Delta^1)$. From the short exact sequence
  \begin{displaymath}
    p^n \widetilde \Z(\Delta^1) \rightarrow R \rightarrow \Z/p^n
  \end{displaymath}
  it follows that we have a weak equivalence $R \overset{\simeq}
  \rightarrow \Z/p^n$. The normalized chain complex of $R$ has a
  generator $1$ in degree zero and a generator $t$ in degree $1$. The
  differential takes $t$ to $p^n \cdot 1$. The
  normalized chain complex $C_*(R)$ of 
  $\HH(R)$ has a generator of the form $1 \otimes t^{\otimes k}$ in
  degree $2k$ 
  and a generator of the form $t^{\otimes k}$ in degree $2k-1$. The
  Hochschild boundary 
  $b$ takes $t^{\otimes k}$ to $p^n(1 \otimes t^{\otimes k-1})$. It
  follows that 
  $\pi_{2k} \HH(R) = \Z/p^n$ and that the odd homotopy groups of
  $\HH(R)$ are zero. In order to compute cyclic homology of $R$, we
  need to evaluate the Connes boundary operator $B$ on the chains of
  the normalized chain complex of $\HH(R)$. The result is that
  $B(t^{\otimes k}) 
  = k(1 \otimes t^{\otimes k})$, and that $B(1 \otimes t^{\otimes k})
  = 0$. It is not easy 
  to compute the higher homology of the bicomplex $(B(R),b,B)$ with
  $B(R)_{s,t} = C_{t-s}$ and with vertical and horizontal differential
  induced by $b$ and $B$ respectively. In degrees up to $2p-1$ the
  horizontal nonzero differentials become isomorphisms after tensoring with
  $\Z/p$. Therefore we have that
  the homology of the total complex of $B(R) \otimes \Z/p$ is a copy
  of $\Z/p$ in degree $k$ when $0 \le k \le 2p-1$. We can
  conclude that if $0 \le i \le p-1$ then $HC_{2i}(R)$ is a cyclic
  $p$-group and $HC_{2i+1}(R) = 0$. To find the order of $HC_{2i}(R)$
  we can consider the spectral sequence associated to the bicomplex
  $(B(R),b,B)$ with $E^1$-term $E^1_{s,t} = \HH_{t-s}(R)$. This
  spectral sequence is concentrated in even total degrees, and
  therefore there are no nonzero differentials. We know that in
  degrees up to $2p-2$ all extensions are maximally nontrivial, and we
  can read off the stated value of $HC_{i}(R)$.

  To see that the map $\widetilde{HC}_i(\Z/p^n) \rightarrow
  \widetilde{HC}_i (\Z/p^{n-1})$ is onto when $0 \le i \le 2p-1$ it
  suffices to check that generators for the group $\Z/p^{(n-1)i} \cong
  \widetilde {HC}_{2i}(\Z/p^{n-1})$ are in the image. This is easy to
  see from the induced map of $E^{\infty}$-terms of the spectral
  sequence considered above.
\end{proof}
\begin{lem}
  The map $\pi_i TC(\Z/p^{n},p)
  \rightarrow 
  \pi_i TC(\Z/p^{n-1},p)$ is onto for $1 \le i \le p-3$ and $n \ge
  2$. Furthermore $\pi_{2j} TC(\Z/p^n,p) = 0$ 
  for $2 \le 2j \le p-3$. 
\end{lem}
\begin{proof}
  The proof goes by induction on $n$. Suppose that $\pi_{2j}
  TC(\Z/p^{n-1},p) = 0$ for $ 2 \le 2j \le p-3$. (By the computation of
  $\pi_* TC(\Z/p,p)$ in \cite[thm. B]{Hesselholt-Madsen} this is true
  for $n = 
  1$.) We have a cofibration 
  sequence
  \begin{displaymath}
    TC(\Z/p^{n},p^{n-1} \Z/p^{n},p) \rightarrow TC(\Z/p^{n},p)
    \rightarrow TC(\Z/p^{n-1},p).
  \end{displaymath}
  Applying proposition \ref{derivedcycliccom}, theorem \ref{mainthm}
  with $I = p^{n-1} \Z/p^n$ and $m=2$ and theorem \ref{relativetc} we
  find that  
  \begin{displaymath}
    \pi_i TC(\Z/p^{n}, p^{n-1} \Z/p^{n},p) \cong
    \begin{cases}
      \Z/p^j & \text{when $i = 2j-1 \le p-3$} \\
      0      & \text{when $i \le p-3$ is even.}
    \end{cases}
  \end{displaymath}
  The statement of the lemma can be read off from the long exact
  sequence associated to the cofibration sequence.
\end{proof}

\begin{cor}
\label{kofzpsquare}
  For $1 \le i \le p-3$, the $K$-groups of $\Z/p^n$ are:
  \begin{displaymath}
     \pi_i K(\Z/p^n) \cong
    \begin{cases}
      0 & \text{if $i$ is even} \\
      \Z/p^{j(n-1)} (p^j -1) & \text{if $i =2j-1$}
    \end{cases}
  \end{displaymath}
\end{cor}

\begin{proof}
  It follows from the above lemma that the map 
  $$\lim_{n} \pi_i
  TC(\Z/p^n,p) \rightarrow \pi_i TC(\Z/p^n,p)$$ is onto and that $\lim_{n}^1
  \pi_i TC(\Z/p^n,p) = 0$ (see \cite{Bousfield-Kan} chap IX and
  XI). We have that 
  $\holim{n} TC(\Z/p^n,p) \simeq TC(\Z^{\wedge}_p,p)$ (see for example
  \cite[thm. 6.1]{Hesselholt-Madsen}), and it follows that the map $\pi_i
  TC(\Z^{\wedge}_p,p) 
  \rightarrow TC(\Z/p^n,p)$ is onto. In \cite{Bokstedt-Madsen} B\"okstedt and
  Madsen have computed 
  $\pi_* TC(\Z^{\wedge}_p,p)$. In the low degrees we are interested in
  it is $\Z^{\wedge}_p$ in odd degrees and $0$ in even strictly
  positive degrees. It follows that the group $\pi_i TC(\Z/p^n,p)$ is
  cyclic. Using the cofibration displayed in the proof of the above
  lemma and the 
  computation of $\pi_* TC(\Z/p,p)$ given in \cite{Hesselholt-Madsen} we
  can by induction prove that 
  \begin{displaymath}
    \pi_i TC(\Z/p^n,p) \cong
    \begin{cases}
      0 & \text{if $0<i \le p-3$ even} \\
      \Z/p^{j(n-1)} & \text{ if $i = 2j-1 \le p-3$.}
    \end{cases}
  \end{displaymath}
  The statement of corollary \ref{kofzpsquare} now follows from McCarthy's
  theorem \ref{relativetc} and from Quillen's computation of $K(\Z/p)$ in
  \cite{QuillenK}. 
\end{proof}

\Addresses\recd

\end{document}